\documentclass[10pt]{article}
\usepackage[T1]{fontenc} 
\usepackage[utf8]{inputenc}

\usepackage{lmodern}

\usepackage{amssymb} 
\usepackage{latexsym}
\usepackage{epsfig}
\usepackage{lscape,graphicx,color,psfrag}
\usepackage{calrsfs}
\usepackage[cmex10]{amsmath}
\usepackage{amsthm}

\theoremstyle{plain}
\newtheorem{lemma}{Lemma}
\newtheorem{assumption}{Assumption}
\newtheorem{proposition}{Proposition}
\newtheorem{theorem}{Theorem}
\newtheorem{corollary}{Corollary}
\theoremstyle{definition}
\newtheorem{definition}{Definition}
\theoremstyle{remark}
\newtheorem{remark}{Remark}

\def\maxim{\mathop{\textup{maximize}}}

\def\st{\mathop{\textup{subject to}}}
\def\diff{\mathop{\Delta}}

\newcommand{\Prob}{\ensuremath{\mathsf{P}}}
\newcommand{\Ex}{\ensuremath{\mathsf{E}}}

\newcommand{\sgn}{\mathop{\mathrm{sgn}}}

\usepackage{appendix}

\usepackage[top=1in, bottom=1.25in, left=0.9in, right=0.9in]{geometry}

\usepackage[sort&compress]{natbib} 
\setcitestyle{authoryear,open={(},close={)}}

\usepackage[colorlinks=true,breaklinks=true,bookmarks=true,urlcolor=blue,
     citecolor=blue,linkcolor=blue,bookmarksopen=false,draft=false]{hyperref}

\def\EMAIL#1{\href{mailto:#1}{#1}}
\def\URL#1{\href{#1}{#1}}         

\begin{document}

\title{A verification theorem for threshold-indexability of real-state  discounted restless bandits}


\author{Jos\'e Ni\~no-Mora\thanks{\small This work has been developed over a number of years since 2008, with support from the Spanish
Ministry of Education and Science [Grant MTM2007-63140], Ministry of Science and Innovation [Grant
MTM2010-20808], and Ministry of Economics and Competitiveness [Grant ECO2015-66593-P].
The author has presented parts of this work in a 2016 seminar at Xerox Research Centre (Grenoble, France) and at the 4th Euro-NGI Conference on Next Generation Internet Networks (NGI 2008) (Krakow, Poland), the 3rd Euro-NF Conference on Network Control and Optimization (NETCOOP 2009) (Eindhoven, Netherlands), the 21st International Symposium on Mathematical
Programming (ISMP 2012) (Berlin, Germany), the 23rd International Symposium on Mathematical Programming (ISMP 2018)
(Bordeaux, France), and the 2013 and 2017 Conferences of the Spanish Royal Mathematical Society. A previous version was
posted on arXiv and submitted in 2015.} 
\\ Department of Statistics \\
    Carlos III University of Madrid \\
     28903 Getafe (Madrid), Spain \\  \EMAIL{jose.nino@uc3m.es}, \URL{http://orcid.org/0000-0002-2172-3983} }
\date{Published in \textit{Mathematics of Operations Research}, vol.\ 45, 465--496, 2020 \\
DOI: \href{https://doi.org/10.1287/moor.2019.0998}{10.1287/moor.2019.0998} \\
Online companion available at \URL{https://pubsonline.informs.org/doi/suppl/10.1287/moor.2019.0998} 
}

\maketitle

\begin{abstract}%
The Whittle index, which characterizes optimal policies for controlling certain
single restless bandit projects (a Markov decision process with two actions: active and
passive) is the basis for a widely used heuristic index policy for the intractable restless
multiarmed bandit problem. Yet two roadblocks need to be overcome to apply such a
policy: the individual projects in the model at hand must be shown to be indexable, so that
they possess a Whittle index; and the index must be evaluated. Such roadblocks can be
especially vexing when project state spaces are real intervals, as in recent sensor scheduling
applications. This paper presents sufficient conditions for indexability (relative to a
generalized Whittle index) of general real-state discrete-time restless bandits under the
discounted criterion, which are not based on elucidating properties of the optimal value
function and do not require proving beforehand optimality of threshold policies as in
prevailing approaches. The main contribution is a verification theorem establishing that, if
project performance metrics under threshold policies and an explicitly defined marginal
productivity (MP) index satisfy three conditions, then the project is indexable with its
generalized Whittle index being given by the MP index, and threshold policies are optimal
for dynamic project control.
\end{abstract}%


\tableofcontents
\newpage
\section{Introduction.}
\label{s:intro}
\subsection{RMABP and Whittle's index policy.}
\label{s:mabpr}
In the \emph{restless multiarmed bandit problem} (RMABP), introduced by \citet{whit88b}, a decision maker aims to maximize the expected reward earned from a finite set of projects.
These are modeled as \emph{restless bandits} (throughout, the terms \emph{project} and \emph{bandit} are used interchangeably), i.e., infinite-horizon \emph{Markov decision processes} (MDPs) ---see \citet{put94}--- that can be operated in two modes, \emph{active} and \emph{passive}.
In the RMABP, a fixed number of projects are selected to be active at each time.
Note that projects are \emph{restless} in that they can change state while passive.

The RMABP extends the classic \emph{MABP}, where projects are \emph{nonrestless} and
one project is engaged at each time, which is a rare example of a multidimensional  MDP for which a tractable policy is available. 
Thus, one can define for each project a scalar function of its state, the   \emph{Gittins index}, so that the resulting \emph{priority-index policy}, which engages at each time a project of largest index, is optimal. See, e.g.,  \citet{gi79}, \citet{whittle80}, and \citet{katehVein87}. \citet{cowanKate15} have extended the scope of such a landmark result to a broader  MABP setting, which allows real-state projects and non-uniform discounting, among other features.
In contrast, index policies are generally suboptimal for the RMABP, as this is known to be PSPACE-hard. See \citet{papTsik99}.

However, \citet{whit88b} extended in part the index-based solution
 approach to the RMABP. He defined an index for restless bandits ---the \emph{Whittle index}--- and proposed as a  heuristic the resulting priority-index policy, which engages at each time the required number of projects with larger indices. He conjectured the asymptotic
optimality of such a policy as both the number of active projects and their total number grow to infinity in fixed ratio, which has been proved by
\citet{wewe90}, \citet{ouyangetal16} and \citet{verloop16} under certain conditions. 
 
Being an intuitive practical heuristic, researchers have deployed Whittle's  policy over the last two decades in diverse applications, providing a growing body of  evidence that its  performance is often nearly optimal. Examples include, e.g., machine maintenance (\citet[Ch.\ 14.6]{whit96}), scheduling multiclass queues (\citet{anselletal03}, \citet{vewe96}, \citet{nmmor06}), 
 admission and routing to queues (\citet{nmmp02,nmcor12,nmcor19}), projects with response delays (\citet{caroYoo10}), dynamic spectrum access (\citet{liuZhao10},  \citet{ouyangetal15}), demand response (\citet{taylorMat14}),  and web crawling (\citet{avraBorkar18}). 

\subsection{Indexability.}
\label{s:indxb} Yet, to apply Whittle's index policy to a given RMABP model researchers need to overcome two roadblocks: (1) establish \emph{indexability} (existence of the index) for 
 the projects in the model at hand;
and (2) devise an efficient index-computing scheme.

Indeed, unlike the Gittins index, which is defined for any classic bandit, the Whittle index is defined only for a limited class of restless bandits, called \emph{indexable}.  \citet[p.\ 292]{whit88b} pointed out that  ``\textit{Indexability cannot be taken for granted $\ldots$ One would very much like to have simple sufficient conditions for indexability$;$ at the moment$,$ none are known}.''
\citet[p.\ 5547]{liuZhao10} further noted that ``\textit{indexability $\ldots$ is often difficult to establish
$\ldots$ and computing Whittle's index can be complex$.$}''

We next outline the 
concepts of indexability and Whittle index as considered herein, which apply to a \emph{single} project.
Consider a project evolving in discrete time that consumes a generic resource. 
At the start of each time period, the controller observes the project state and then chooses
 an action (active or passive).
Then the project
(i)  consumes an amount of  resource and yields   a  
reward, discounted over time;
and (ii) 
moves  to the next state according to a Markov transition
law. 
  
Action choice is governed by a policy $\pi$ from the
class $\Pi$ of \emph{nonanticipative policies}.
Suppose that the resource is charged at a unit price 
$\lambda$, and let $V_\lambda(x, \pi)$ be the resulting net project value starting from state $x$ under $\pi$, giving the expected total discounted value of rewards earned minus resource charges incurred.
Now, consider the parametric collection of \emph{$\lambda$-price problems} $\{P_\lambda\colon \lambda \in \mathbb{R}\}$, with
\begin{equation}
\label{eq:0lpricepfg1}
P_\lambda\colon \qquad \textup{find } \pi_\lambda^* \in \Pi \textup{ such that } V_\lambda(x, \pi_\lambda^*) = V_\lambda^*(x) \triangleq \sup_{\pi \in \Pi} \, V_\lambda(x, \pi)  \textup{ for every state }  x.
\end{equation}
Thus, a \emph{$P_\lambda$-optimal} policy achieves an optimal tradeoff between rewards and resource usage costs. 

We shall call the project \emph{indexable} if there exists a scalar map $\lambda^*(x)$ on $\mathsf{X}$ (the project's \emph{Whittle index}) characterizing $P_\lambda$-optimal policies for any price $\lambda$, as follows: when the project is in state $x$, it is $P_\lambda$-optimal to take
the active (resp.\ passive) action   if and only if 
$\lambda^*(x) \geqslant \lambda$ (resp.\ $\lambda^*(x) \leqslant \lambda$).
 
\subsection{Real-state projects: prevailing approaches to indexability based on optimal thresholds.}
\label{s:rmabprsp}
While most work on the RMABP has considered discrete-state projects,  models with real-state projects have drawn growing attention since the last decade, motivated by \emph{partially observed MDP} (POMDP) formulations of sensor scheduling applications. See, e.g.,  \citet[Ch.\ 8]{krishn16}  and the surveys \citet{washburn08} and \citet{kuhnNaz17}.

In such applications, a set of sensors are to be dynamically allocated to sense a larger set of objects, e.g., moving targets (in multitarget tracking, where the sensors are radars) or communication channels (in multichannel access, where sensors probe channel transmission availability). Due to the limited number of sensors and, possibly, their providing noisy measurements, true object states are unknown at decision times.
Instead, under a POMDP model sensor allocation decisions are based on the \emph{belief state} of each object, which gives the posterior state distribution   and is updated by a Bayesian filter depending on whether the object is sensed.
Belief states can take a continuum of values, being scalar in certain models.
Thus, identifying objects with projects, the active and passive actions with sensing and not sensing an object, and using belief states as project states, when the latter are scalar we obtain special cases of the RMABP with real-state projects.

\citet{lascala06} first formulated the problem of dynamic radar scheduling to minimize the track error variance of two moving targets as an RMABP, where projects follow deterministic Kalman filter dynamics.
They proved optimality of the greedy policy in the symmetric case. 

More recent work has addressed indexability in POMDP models with real-state projects, 
including, e.g., multisite tracking 
(\citet{lenyferon08}), dynamic multichannel access (\citet{liuZhao10} and \citet{ouyangetal15}), continuous-time (\citet{lenyetal11}) and discrete-time (\citet{danceSi15}) multitarget tracking, demand response (\citet{taylorMat14}),  web crawling (\citet{avraBorkar18}), and hidden Markov bandits (\citet{meshrametal18}). 

Such work follows a common pattern, based on  elucidating and exploiting properties of the optimal value function $V_\lambda^*(x)$ of $\lambda$-price problems (\ref{eq:0lpricepfg1}) for the model at hand through \emph{dynamic programming} (DP) arguments, to prove optimality of \emph{threshold policies}. 
Such policies, also known as \emph{control limit policies}, make the project active when its state is (at or) above a given threshold, and make it passive otherwise. They are widely used in applications due to their simplicity and ease of implementation.  

A typical analysis proceeds in two stages.
The first is to identify and exploit properties of  $V_\lambda^*(x)$ to prove existence of an optimal threshold policy for each $\lambda$-price problem, with the optimal threshold $z^*(\lambda)$ being unique for a certain interval of $\lambda$'s.
In the second stage, it is claimed that to establish indexability it suffices to show that $z^*(\lambda)$ is \emph{increasing} (or, in some papers, \emph{nondecreasing}) in $\lambda$. The Whittle index $\lambda^*(x)$ is  then obtained by inverting $z^*(\lambda)$, i.e., solving the equation $z^*(\lambda) = x$ in $\lambda$.

Yet such analyses, while insightful and useful, are often incomplete, due to: (1) as proving optimality of threshold policies may be elusive, it is sometimes only conjectured (as in \citet{lenyetal11} and \citet{danceSi15}), or the arguments provided are sketchy and not fully rigorous; and
(2) showing \emph{increasingness} (in the strict sense) of $z^*(\lambda)$ does not suffice to establish indexability; \emph{continuity} of $z^*(\lambda)$ needs also be proven (see Proposition \ref{pro:iwrttpndlc}).

Indeed, proving optimality of threshold policies can be far from trivial. 
It can be formulated in terms of $V_\lambda^{\langle a, *\rangle}(x)$, the optimal value with initial action $a$. 
Thus, it needs to be shown that, for $\lambda$ in a certain interval, the equation $\Delta_{a=0}^{a=1} V_\lambda^{\langle a, *\rangle}(x) \triangleq V_\lambda^{\langle 1, *\rangle}(x) - V_\lambda^{\langle 0, *\rangle}(x)  = 0$ in state variable $x$ has a unique root $z^*(\lambda)$, with 
$\Delta_{a=0}^{a=1} V_\lambda^{\langle a, *\rangle}(x) < 0$ for $x < z^*(\lambda)$ and $\Delta_{a=0}^{a=1} V_\lambda^{\langle a, *\rangle}(x) > 0$ for $x > z^*(\lambda)$.

Such a result follows if the difference $V_\lambda^{\langle 1, *\rangle}(x) - V_\lambda^{\langle 0, *\rangle}(x)$ is concave in $x$ and changes sign at the endpoints of the state space, as in \citet{lenyferon08} and   
\citet{liuZhao10}, where $V_\lambda^{\langle 1, *\rangle}(x)$ is affine and $V_\lambda^{\langle 0, *\rangle}(x)$ convex.
Yet, in the model extension in 
 \citet{ouyangetal15}, while  the $V_\lambda^{\langle a, *\rangle}(x)$ are convex in $x$, it is no longer ensured that $\Delta_{a=0}^{a=1} V_\lambda^{\langle a, *\rangle}(x)$ is concave, as $V_\lambda^{\langle 1, *\rangle}(x)$ may not be affine.  Hence,  existence of a $z^*(\lambda)$ as above 
cannot be concluded from $\Delta_{a=0}^{a=1} V_\lambda^{\langle a, *\rangle}(x)$ changing sign.
 
The most widely used approach for proving optimality of threshold policies draws on 
the theory of \emph{supermodularity}. See, e.g., \citet{topkis98}, \citet[Ch.\ 8.2]{heymSob84ii}, \citet[Ch.\ 6.11.2]{put94}, \citet{altmanetal95}, and \citet[Ch.\ 9]{krishn16}.
Thus, if $V_\lambda^{\langle a, *\rangle}(x)$ is \emph{strictly supermodular} in $(x, a)$, so $\Delta_{a=0}^{a=1} V_\lambda^{\langle a, *\rangle}(x)$ is increasing in $x$, being negative for small $x$  and positive for  large $x$, existence of a unique threshold $z^*(\lambda)$ satisfying the required properties follows.

Yet, in the setting of a Kalman filter target tracking model, \citet{danceSi19} identify an instance demonstrating that efforts to prove optimality of threshold policies in such a model through the above approaches are doomed to failure. Using the notation in Section \ref{s:mdf}, 
their instance has rewards $r(x, a) = - x$, passive dynamics $X_{t+1} = X_t + 1$,  active dynamics 
$X_{t+1} = 1/(\alpha + 1/(X_t+1))$ with $\alpha = 0.1$, and discount factor $\beta = 0.95$. Figure 1 in their paper shows that, for $\lambda = 0.7647$, $\Delta_{a=0}^{a=1} V_\lambda^{\langle a, *\rangle}(x)$ is neither concave nor monotonic  over  $x \in [0, 0.25]$.
The author has independently verified such a result, which highlights a limitation of the scope of prevailing approaches to optimality of threshold policies and indexability, raising the need for different  approaches with a wider scope.

\subsection{PCL-based sufficient indexability conditions.}
\label{s:pclbsicrsp} 
An alternative to prevailing approaches to indexability for discrete-state restless bandits was introduced and developed in \citet{nmaap01,nmmp02,nmmor06,nmtop07} based on two
sufficient conditions, which can be verified analytically or numerically, along with an index algorithm.
The validity of such conditions stems from satisfaction by project metrics of  \emph{partial conservation laws} (PCLs),  
whence we shall refer to \emph{PCL-indexability conditions}. PCLs extend earlier types of  conservation laws in stochastic scheduling, which have been used to prove  optimality of index policies through achievable performance region analyses. See, 
e.g., \citet{coffMitrani80}, \citet{shanthiYao92}, \citet{nm96}, and \citet{nmconsLawsEORMS}.  
Yet, the proof techniques used in the aforementioned work rely critically on the discreteness of the project state space, so their extension to real-state projects is nontrivial. 

\citet{nmngi08} first outlined, without proving its validity, an extension to real-state projects of such PCL-indexability conditions, which roughly correspond  to (PCLI1) and (PCLI2) in this paper (see Definition \ref{def:pcli}).  Yet condition (PCLI3) herein was not considered, nor was the continuity requirement in (PCLI2).
That paper further illustrated application of the conditions in the model in \citet{lenyferon08} and \citet{liuZhao10}, sketching an incomplete analysis.

Later work exploring use of such an approach includes
  \citet{nmnetcoop09,nm16} and the PhD thesis of \citet{vill12}. The latter, supervised by the author, deepened the study of the models in
  \citet{nmsv09,nmv11}, including the multitarget Kalman filter tracking model in \citet{lascala06}.  
While such work did not establish validity of the conditions, it presented partial analyses as well as numerical evidence supporting the practical value of the proposed approach: the Whittle index could be efficiently computed and the index policy was found to be nearly optimal. 

In the first arXiv version of this paper, \citet{nmvt15}
modified the PCL-indexability conditions in \citet{nmngi08} by adding a continuity requirement in (PCLI2), as well as the new condition  (PCLI3) in Definition \ref{def:pcli}. \citet{nmvt15}  gave a  proof that conditions (PCLI1--PCLI3) imply indexability \emph{and} optimality of threshold policies for the $\lambda$-price problems. The proof relied on analysis of an infinite-dimensional 
\emph{linear programming} (LP) formulation of such problems. 

Based on \citet{nmvt15}, \citet{danceSi19} have applied 
the PCL-indexability conditions herein to establish for the first time indexability, along with optimality of threshold policies, of an extension of the aforementioned Kalman filter target tracking model.
To prove (PCLI1--PCLI3) they deploy a novel analysis of action trajectories under threshold policies drawing on the theory of combinatorics on words. In a preliminary arXiv version of such work, 
 further outline a short proof of the validity of such conditions for their model, which motivated the author to develop the simpler proof presented here for the general case that avoids the infinite-dimensional LP analysis in  \citet{nmvt15}.

\subsection{Goal and main results.}
\label{s:gmr} 
The major goal of this paper is to develop a rigorous PCL-based approach to establish indexability and evaluate the Whittle index of real-state discounted restless bandits in a form consistent with optimality of threshold policies for their $\lambda$-price problems. 

The main result is a \emph{verification theorem} (Theorem \ref{the:pcliii})  which, under three \emph{PCL-indexability conditions},  ensures in one fell swoop  \emph{both} indexability of a real-state restless bandit \emph{and} optimality of threshold policies for its $\lambda$-price problems (\ref{eq:0lpricepfg1}).
This differs from prevailing approaches (see Section \ref{s:rmabprsp}), where first optimality of threshold policies is proved, and then indexability is   established. 

Under the proposed conditions, 
 the Whittle index is given by a \emph{marginal productivity} (MP) index that is explicitly defined in terms of  project metrics under threshold policies. 
This is in contrast both to the Whittle index definition, which 
 characterizes it implicitly, and to 
 the way in which the index is determined in prevailing approaches, by inverting an optimal threshold function.
 
The conditions are formulated as properties
 of \emph{project metrics} under threshold policies.
Thus, we consider the \emph{reward} and \emph{resource $($usage$)$} metrics denoted respectively by $F(x, z)$ and $G(x, z)$, with $x$ the initial state and $z$  shorthand for the \emph{$z$-policy} (see Section \ref{s:assub}).
We also consider
\emph{marginal reward} and \emph{marginal resource} metrics, denoted respectively by $f(x, z) \triangleq F(x, \langle 1, z\rangle) - F(x, \langle 0, z\rangle)$ and $g(x, z) \triangleq G(x, \langle 1, z\rangle) - G(x, \langle 0, z\rangle)$, where 
policy $\langle a, z\rangle$ first takes action $a$ and then follows the $z$-policy.
We further consider the \emph{marginal productivity} (MP) metric $m(x, z) \triangleq f(x, z)/g(x, z)$, which is defined when $g(x, z) \neq 0$, and the \emph{MP index}, given by $m(x) \triangleq m(x, x)$ for every state $x$.

The \emph{PCL-indexability conditions} are: (PCLI1) the marginal resource metric $g(x, z)$ is positive for any state $x$ and threshold $z$;
(PCLI2) the MP index $m(x)$ is monotone nondecreasing and continuous; and 
(PCLI3) the reward and resource metrics and the MP index are linked by
\begin{equation}
\label{eq:pcli3}
F(x, z_2) -  F(x, z_1) = \int_{(z_1, z_2]} m(z) \, G(x, dz), \quad -\infty < z_1 < z_2   < \infty,
\end{equation}
so, for each $x$, and viewed as functions of the threshold $z$, $F(x, z)$ is an  indefinite \emph{Lebesgue--Stieltjes} (LS) \emph{integral} (see, e.g.,  \citet{cartvanBrunt00}) of the MP index $m(z)$ with respect to $G(x, z)$.
Proposition \ref{pro:mpirnd}  shows that (PCLI3) characterizes  the index as a \emph{Radon--Nikodym derivative}.

Theorem \ref{the:pcliii} ensures that a project satisfying conditions (PCLI1--PCLI3) is indexable with $m(x)$ being its Whittle index. 
Furthermore, as the latter is nondecreasing and continuous, the optimality of threshold policies for the project's $\lambda$-price problems (\ref{eq:0lpricepfg1}) follows. See Lemma \ref{lma:iwrttpndlc}.

This result reduces the task of establishing a project's indexability to that of proving certain properties of project metrics \emph{under threshold policies}, and  provides an explicit expression for the index.
Such an approach is in contrast to prevailing approaches, based on exploiting properties of the optimal value function, which require establishing beforehand optimality of threshold policies. 

Lemma \ref{lma:iwrttpndlc} and Proposition \ref{lma:intgradz} show that conditions (PCLI2, PCLI3) are necessary  for \emph{threshold-indexability} (see Definition \ref{def:indxtp}), and that the Whittle index matches the MP index under (PCLI1). See Proposition \ref{lma:phiwphimpi}(c).
Proposition \ref{pro:cfhMPI} gives an index-computing recursion.

Note that the PCL-indexability conditions for discrete-state projects referred to in Section \ref{s:pclbsicrsp} roughly correspond to (PCLI1, PCLI2), as in such a setting they imply (PCLI3).
In the present real-state setting, since (PCLI3) may be hard to verify, 
we give Propositions \ref{pro:pwdiffFG} and \ref{pro:1cpcli3}, which ensure that (PCLI3) holds under (PCLI1, PCLI2) and some additional simpler conditions.

\subsection{Organization of the paper.}
\label{s:porg}
The remainder of the paper is organized as follows.
Section \ref{s:mdf} sets up notation and 
presents concepts  to formulate the main result in Section \ref{s:pclbsic}.
Section \ref{s:potip} motivates the PCL-indexability conditions.
Section \ref{s:srsti} presents results for the indexability analysis of \emph{strongly thresholdable} projects (see Definition \ref{def:spotp}), which are relevant for complementing the approaches reviewed in Section \ref{s:rmabprsp}.
Section \ref{s:saom} reformulates the $\lambda$-price problem 
$P_\lambda$ as an LP problem over a space of measures.
Section \ref{s:sdr} presents a key result on decomposition of project  metrics.
Section \ref{s:ppmftv} establishes required properties of project metrics as functions 
of the threshold, including that they are \emph{c\`adl\`ag} (``continue \`a droite, limite \`a gauche,'' i.e., right-continuous with left limits), and characterizes the MP index as a Radon--Nikodym derivative. 
Section \ref{s:opicfpmpm} derives further required relations between project metrics. 
Section \ref{s:ptpcliii} gives the proof of Theorem \ref{the:pcliii}.
Section \ref{s:ccpcli3} presents tools to establish that condition
 (PCLI3) holds.
 Section \ref{s:pmaic} addresses index and metrics computation.
Section \ref{s:ex} illustrates application of the proposed approach in two 
models.
Section \ref{s:conc} concludes.

An online companion containing two appendices provides relevant ancillary material. Section A motivates
the interest of indexability through a general resource-constrained RMABP. Section B presents geometric and
economic interpretations of the MP index under PCL-indexability, showing in particular that it can be
characterized as a \emph{resource shadow price}.

\section{Preliminaries and formulation of the main result.}
\label{s:mdf}
We next set up notation,  
describe the single-project model, and present required concepts to formulate the main result. 

\subsection{Notation.}
\label{s:not}
We denote by $\mathsf{X} \triangleq \{x \in \mathbb{R}\colon \ell \leqslant x \leqslant  u\}$ the project state space, which is a closed (possibly unbounded) interval of $\mathbb{R}$ with endpoints $-\infty \leqslant \ell < u \leqslant \infty$, and by $\mathsf{K}  \triangleq \mathsf{X} \times \{0, 1\}$ its state-action space.
For $S \in  \mathcal{B}(\mathsf{X})$, the Borel $\sigma$-algebra of $\mathsf{X}$, we write the complement of $S$ in $\mathsf{X}$ as $S^c \triangleq
\mathsf{X} \setminus S$, and denote 
	by $1_S(\cdot)$ its indicator function. 
All subsets of $\mathbb{R}$ and functions referred to are assumed to be Borel measurable.
We write functions as, e.g., $J(\cdot)$, and denote by $|J|(\cdot)$ the corresponding absolute value, so $|J|(x) = |J(x)|$.
We denote by $\mathbb{C}(\mathsf{X})$ the space of real-valued continuous functions on $\mathsf{X}$.

Throughout, $x$ and $y$ denote project states, and 
we use $z \in \overline{\mathbb{R}} \triangleq \mathbb{R} \cup \{-\infty, \infty\}$ to denote thresholds. For a given threshold $z$, we use shorthands such as, e.g., $\{x > z\}$ to mean $\{x \in \mathsf{X}\colon x > z\}$.

We write as $\int_B J \, d\alpha$ or $\int_B J(x) \, \alpha(dx)$ the LS integral (see \citet{cartvanBrunt00}) of a function $J(\cdot)$ with respect to a function $\alpha(\cdot)$ over  
$B$. When $B$ is the domain of $J(\cdot)$ we write $\int J \, d\alpha$.

We use the terms ``increasing'' and ``decreasing'' in the strict sense, and ``iff'' to mean ``if and only if.'' 

\subsection{Project model,  performance metrics, and  $\lambda$-price problems.}
\label{s:rrupm}
We consider a general restless bandit model of optimal dynamic resource allocation to a  
project, whose state $X_{t}$ moves  in discrete time over an infinite  horizon across the state space $\mathsf{X}$. 

At the start of period $t = 0, 1, \ldots$, a controller observes the project state $X_t$ and then chooses
 \emph{action} $A_t \in \{0, 1\}$,
where $A_t = 1$ (resp.\ $A_t = 0$) codes the project being  \emph{active} (resp.\  \emph{passive}). When $(X_t, A_t) = (x, a)$: 
(i) the project consumes 
$c(x, a)$  units of a generic resource and yields   a  
reward $r(x, a)$, discounted with factor $\beta \in [0, 1)$;
and (ii) 
its  state moves  to $X_{t+1}$  according to a Markov transition
kernel $\kappa^a(\cdot, \cdot)$, so 
$\Prob\{X_{t+1} \in B \, | \, (X_t, A_t) = (x, a)\} =
 \kappa^a(x, B)$ for $B \in \mathcal{B}(\mathsf{X})$.

We make the following standing assumption. 
\begin{assumption}
\label{ass:first} For each action $a \in \{0, 1\}$\textup{:}
\begin{itemize}
\item[\textup{(i)}] $r(\cdot, a), c(\cdot, a) \in \mathbb{C}(\mathsf{X}),$ with $c(\cdot, a)$ satisfying
$0 \leqslant c(x, 0) < c(x, 1)$ for $x \in \mathsf{X};$ 
\item[\textup{(ii)}] there is  a function 
   $w\colon \mathsf{X} \to [1, \infty)$ and $M > 0, \beta \leqslant \gamma < 1$
 such that$,$ for any state $x,$
\begin{itemize}
\item[\textup{(ii.a)}] $\max \{\vert r \vert(x, a),  c(x, a)\} \leqslant M \, w(x);$ and
\item[\textup{(ii.b)}] $\beta \widetilde{w}(x, a) \leqslant \gamma
  \, w(x),$ where $\widetilde{w}(x, a) \triangleq \int w(y) \, \kappa^a(x, dy);$
\end{itemize}
\end{itemize}
\end{assumption}

\begin{remark}
\label{re:first} 
\begin{itemize}
\item[\textup{(i)}] Assumption \ref{ass:first}(ii) corresponds to the  \emph{weighted sup-norm} framework for discounted MDPs (see, e.g.,  
\citet{lippman75,wessels77} and \citet[assumption 8.3.2]{herlerLass99}), which allows functions $r(\cdot, a)$ and $c(\cdot, a)$ to be unbounded. 
\item[\textup{(ii)}] If $\widetilde{w}(x, a) \leqslant w(x)$, Assumption \ref{ass:first}(ii.b) holds for any $\beta$ by taking $\gamma \triangleq \beta$.
\end{itemize}
\end{remark}

Actions are selected through 
a \emph{policy} $\pi = (\pi_t)_{t=0}^\infty$ from the
class $\Pi$ of \emph{nonanticipative policies}, which at time $t$ selects action $A_t = a$ with probability 
$\pi_t(a \, | \, \mathcal{H}_t)$, a measurable function of the history $\mathcal{H}_t \triangleq \big((X_s, A_s)_{s=0}^{t-1}; X_t\big)$.
We shall further refer to the class 
$\Pi^{\scriptscriptstyle \textup{SR}}$ of \emph{stationary randomized policies}, where $\pi_t(a \, | \, \mathcal{H}_t) = q(a \, | \, X_t)$ for some $q(\cdot \, | \, \cdot) \geqslant 0$ with $\sum_{a} q(a \, | \, x) = 1$, 
and to the class $\Pi^{\scriptscriptstyle \textup{SD}}$ of \emph{stationary deterministic policies}, where $A_t = 1_B(X_t)$ for some \emph{active region} $B \in \mathcal{B}(\mathsf{X})$, in which  
case we shall refer to the \emph{$B$-policy}, writing $\pi = B$.

For a policy  $\pi \in \Pi$ and an initial-state probability measure $\nu_0$ on $\mathcal{B}(\mathsf{X})$  (written $X_0 \sim \nu_0$), let 
$\Prob_{\nu_0}^{\pi}$ be the probability measure determined (by Ionescu-Tulcea's extension theorem; see \citet[proposition C.10]{herlerLass96}) on the space $\mathsf{K}^\infty$ of state-action paths with the product $\sigma$-algebra, with expectation   
$\Ex_{\nu_0}^{\pi}$. We write $\Prob_x^{\pi}$ and $\Ex_x^{\pi}$
when $\nu_0 = \delta_x$, the Dirac measure at $x$.

We evaluate a policy $\pi$ starting from an initial state drawn from $\nu_0$ by the \emph{reward} and  \emph{resource $($\!usage$)$ performance metrics} 
\begin{equation}
\label{eq:FGboldx0pis}
F(\nu_0, \pi) \triangleq \Ex_{\nu_0}^{\pi}\left[\sum_{t = 0}^\infty 
\beta^t r(X_{t}, A_{t})\right] \quad \textup{and} \quad G(\nu_0, \pi)\triangleq \Ex_{\nu_0}^{\pi}\Bigg[\sum_{t = 0}^\infty 
\beta^t c(X_{t}, A_{t})\Bigg],
\end{equation}
which we write as $F(x, \pi)$ and $G(x, \pi)$ when $\nu_0 = \delta_x$, and as  $F(\nu_0, B)$ and $G(\nu_0, B)$
when $\pi = B$.

If resource usage is charged at a price $\lambda$ per unit, denote by $V_\lambda(\nu_0, \pi) \triangleq F(\nu_0, \pi)- \lambda G(\nu_0, \pi)$ the \emph{$($net$)$ project value}, and by  $V_\lambda^*(\nu_0) \triangleq \sup_{\pi \in \Pi} \, V_\lambda(\nu_0, \pi)$ the \emph{optimal project value}.

For any resource price $\lambda$ consider the  
\emph{$\lambda$-price problem} 
\begin{equation}
\label{eq:lpricepfg1}
P_\lambda\colon \qquad \textup{find } \, \pi_\lambda^* \in \Pi \, \textup{ such that } \, V_\lambda(x, \pi_\lambda^*) = V_\lambda^*(x)  \, \textup{ for every }  x \in \mathsf{X}.
\end{equation} 
We call a policy $\pi_\lambda^*$ solving (\ref{eq:lpricepfg1}) \emph{$P_\lambda$-optimal}.
Consider
 the parametric problem collection 
\begin{equation}
\label{eq:ppc}
\mathcal{P} \triangleq \{P_\lambda\colon \lambda \in \mathbb{R}\}.
\end{equation}

We next  review results from the weighted sup-norm framework that are needed in the sequel.
Let $\mathbb{B}_w(\mathsf{K})$ and $\mathbb{B}_w(\mathsf{X})$ be the  Banach spaces (see \citet[proposition 7.2.1]{herlerLass99}) of
 \emph{$w$-bounded} functions $u\colon \mathsf{K} \to \mathbb{R}$ and $v\colon \mathsf{X}  \to \mathbb{R},$ respectively, having finite
\emph{$w$-norms} 
\begin{equation}
\label{eq:vxwnorm}
\| u \|_w \triangleq \sup_{(x, a) \in \mathsf{K}} \, \frac{\vert u \vert (x, a) }{w(x)}
  \quad \textup{and} \quad
\| v \|_w \triangleq \sup_{x \in \mathsf{X}} \frac{\vert v \vert (x)}{w(x)},
\end{equation}
and write  as $\mathbb{P}_w(\mathsf{X})$ the space of probability measures $\nu_0$ on $\mathcal{B}(\mathsf{X})$ with finite 
$w$-norm
\begin{equation}
\label{eq:pmwnorm}
\| \nu_0 \|_w \triangleq \int w \, d\nu_0 = \Ex_{\nu_0}[w(X_0)].
\end{equation}
Let further $\mathbb{C}_w(\mathsf{K}) \triangleq \mathbb{C}(\mathsf{K}) \cap \mathbb{B}_w(\mathsf{K})$ and $\mathbb{C}_w(\mathsf{X}) \triangleq \mathbb{C}(\mathsf{X}) \cap \mathbb{B}_w(\mathsf{X})$ denote the corresponding subspaces of $w$-bounded continuous functions on $\mathsf{K}$ and $\mathsf{X}$.

Henceforth, $F(\cdot, \pi)$, $G(\cdot, \pi)$ and $V_\lambda(\cdot, \pi)$ denote the reward, resource and net value project metrics viewed as functions of the initial state $x$ for a fixed policy $\pi$.

\begin{remark}
\label{re:assboundfs} 
Under Assumption \textup{\ref{ass:first} (}cf.\ \textup{\citet[theorem \textup{8.3.6}]{herlerLass99}):}
\begin{itemize}
\item[\textup{(i)}] For any action $a$, $r(\cdot, a), c(\cdot, a) \in \mathbb{C}_w(\mathsf{X});$  $r(\cdot, \cdot), c(\cdot, \cdot) \in \mathbb{C}_w(\mathsf{K});$ and $\kappa^a(x, \cdot) \in \mathbb{P}_w(\mathsf{X}).$
\item[\textup{(ii)}] For any $\pi \in \Pi$ and $\lambda \in \mathbb{R},$ $F(\cdot, \pi), G(\cdot, \pi), V_\lambda(\cdot, \pi) \in \mathbb{B}_w(\mathsf{X});$ and$,$ with $M_\gamma \triangleq M/(1-\gamma),$ 
\begin{equation}
\label{eq:FGMwnon}
\max \{\| F(\cdot, \pi) \|_w, 
\| G(\cdot, \pi) \|_w\} \leqslant M_\gamma \quad \textup{and} \quad
\| V_\lambda(\cdot, \pi)\|_w \leqslant (1 + |\lambda|) M_\gamma.
\end{equation}
\item[\textup{(iii)}] For $\pi \in \Pi,$ $\nu_0 \in  \mathbb{P}_w(\mathsf{X})$ and $\lambda \in \mathbb{R},$ 
$F(\nu_0, \pi),$ $G(\nu_0, \pi)$ and $V_\lambda(\nu_0, \pi)$ are well defined and finite$,$ being equal to
$\int F(x, \pi) \, \nu_0(dx),$ $\int G(x, \pi) \, \nu_0(dx)$ and $\int V_\lambda(x, \pi) \, \nu_0(dx),$ respectively$,$  and satisfy
\begin{equation}
\label{eq:FGMnuwnon}
\max \{|F|(\nu_0, \pi) , 
G(\nu_0, \pi)\} \leqslant M_\gamma \| \nu_0 \|_w \quad \textup{and} \quad
|V_\lambda|(\nu_0, \pi) \leqslant (1 + |\lambda|) M_\gamma \| \nu_0 \|_w.
\end{equation}
\item[\textup{(iv)}] For any active region $B \in \mathcal{B}(\mathsf{X})$, 
$F(\cdot, B)$ and $G(\cdot, B)$  are  characterized $($cf.\ \citet[remark 8.3.10]{herlerLass99}$)$ as the
unique solutions in $\mathbb{B}_w(\mathsf{X})$ to the fixed-point
equations
\begin{equation}
\label{eq:fxbee}
\begin{split}
F(x, B) & = 
r(x, 1_{B}(x)) + \beta \int F(y, B) \,
  \kappa^{1_B(x)}(x, dy),  \quad x \in \mathsf{X}, \\
G(x, B) & = 
c(x, 1_B(x)) + \beta \int G(y, B) \,
  \kappa^{1_B(x)}(x, dy),  \quad x \in \mathsf{X}.
\end{split}
\end{equation}
Later (see Section \ref{s:sdr}) we will find it convenient to use the reformulation of $(\ref{eq:fxbee})$ 
in terms of the bounded linear operator $\mathcal{L}^\star\colon \mathbb{B}_w(\mathsf{X}) \to \mathbb{B}_w(\mathsf{K})$ mapping $v(\cdot)$  to 
\begin{equation}
\label{eq:tbetastu}
\mathcal{L}^\star v(x, a) \triangleq v(x) - \beta \int v(y) \, 
\kappa^a(x, dy), \quad (x, a) \in \mathsf{K},
\end{equation}
namely$,$ 
\begin{equation}
\label{eq:fseqs}
\begin{split}
\mathcal{L}^\star F(\cdot, B)(x, 1_B(x))  & =  r(x, 1_B(x)),  \quad x \in \mathsf{X},  \\
\mathcal{L}^\star G(\cdot, B)(x, 1_B(x))  & =  c(x, 1_B(x)),  \quad x \in \mathsf{X}.
\end{split}
\end{equation}
\item[\textup{(v)}] For any price $\lambda \in \mathbb{R},$  the \emph{Bellman operator} $\mathcal{T}_\lambda\colon \mathbb{B}_w(\mathsf{X}) \to \mathbb{B}_w(\mathsf{X})$ defined by  
\begin{equation}
\label{eq:Vstast}
\mathcal{T}_\lambda v(x)  \triangleq \max_{a \in \{ 0, 1\}} \, \mathcal{T}_\lambda^a v(x), \quad \textup{where } \mathcal{T}_\lambda^a v(x) \triangleq r(x, a) - \lambda c(x, a) + \beta 
  \int v(y) \, \kappa^a(x, dy), 
  \end{equation} is a \emph{contraction mapping} with modulus $\gamma$ (see Assumption \textup{\ref{ass:first}(ii)}) and the \emph{Bellman equation} (BE) 
\begin{equation}
\label{eq:dpoes}
v = \mathcal{T}_\lambda v,
\end{equation}  
has the \emph{optimal value function} $V_\lambda^*(\cdot)$ as its unique fixed point in $\mathbb{B}_w(\mathsf{X})$, which satisfies 
\begin{equation}
\label{eq:VMwnon}
\| V_\lambda^*\|_w \leqslant (1 + |\lambda|) M_\gamma.
\end{equation}
\item[\textup{(vi)}] For any price $\lambda  \in \mathbb{R}$  there is a $P_\lambda$-optimal stationary deterministic policy; a policy in $\Pi^{\scriptscriptstyle \textup{SD}}$ is $P_\lambda$-optimal iff it selects in each state $x$ an action $a$ with
$\mathcal{T}_\lambda^a V_\lambda^*(x) \geqslant \mathcal{T}_\lambda^{1-a} V_\lambda^*(x)$.
\end{itemize}
\end{remark}

We  denote the optimal net project value starting from $X_0 = x$ with initial action $A_0 = a$ by  
\begin{equation}
\label{eq:vlastarx}
V_\lambda^{\langle a, *\rangle}(x) \triangleq \mathcal{T}_\lambda^a V_\lambda^*(x) =  
r(x, a) - \lambda c(x, a) + \beta
 \int V_\lambda^*(y) \, \kappa^a(x, dy).
\end{equation}

We say that \emph{action $a$ is $P_\lambda$-optimal in state $x$}  if  $V_\lambda^{\langle a, *\rangle}(x) \geqslant V_\lambda^{\langle 1-a, *\rangle}(x)$.
It is convenient to reformulate such a definition in terms of the sign of $\Delta_{a=0}^{a=1} V_\lambda^{\langle a, *\rangle}(x) \triangleq V_\lambda^{\langle 1, *\rangle}(x) - V_\lambda^{\langle 0, *\rangle}(x)$,  the \emph{marginal value of engaging the project in state $x$}.  Thus,  
action $a = 1$ is $P_\lambda$-optimal in $x$ if $\Delta_{a=0}^{a=1} V_\lambda^{\langle a, *\rangle}(x) \geqslant 0$; action $a = 0$ is $P_\lambda$-optimal in $x$ if $\Delta_{a=0}^{a=1} V_\lambda^{\langle a, *\rangle}(x) \leqslant 0$; and both actions are $P_\lambda$-optimal in $x$ if \begin{equation}
\label{eq:bothactopt}
\Delta_{a=0}^{a=1} V_\lambda^{\langle a, *\rangle}(x) = 0.
\end{equation}

\subsection{Indexability.}
\label{s:piip}
We will address the parametric problem collection $\mathcal{P}$ in (\ref{eq:lpricepfg1}) through
 the concept of
\emph{indexability}, extended by \citet{whit88b} from classic to restless bandits with resource consumption  $c(x, a) \triangleq a$, and further extended by \citet{nmmp02} to general $c(x, a)$. 

For a given resource price $\lambda$, the \emph{$P_\lambda$-optimal active} and \emph{passive regions} are given by
\begin{equation}
\label{eq:Sstarlam}
S^{*, 1}_\lambda \triangleq \big\{x \in \mathsf{X}\colon \Delta_{a=0}^{a=1} V_\lambda^{\langle a, *\rangle}(x) \geqslant 0\big\} \quad \textup{and} \quad S^{*, 0}_\lambda \triangleq \big\{x \in \mathsf{X}\colon  \Delta_{a=0}^{a=1} V_\lambda^{\langle a, *\rangle}(x) \leqslant 0\big\},
\end{equation}
respectively. 
Under indexability, such regions are characterized by an 
\emph{index} attached to project states. Note that the definition below refers to the \emph{sign} function $\sgn\colon \mathbb{R} \to \{-1, 0, 1\}$.

\begin{definition}[Indexability]
\label{def:indx}
The project is \emph{indexable} if there is a map $\lambda^*\colon \mathsf{X} \to \mathbb{R}$ with 
\begin{equation}
\label{eq:indxsigndef}
\sgn \Delta_{a=0}^{a=1} V_\lambda^{\langle a, *\rangle}(x) = \sgn (\lambda^*(x) - \lambda), \quad \textup{for } x \in \mathsf{X} \textup{  and  } \lambda \in \mathbb{R}.
\end{equation}
We call $\lambda^*(\cdot)$ the project's \emph{Whittle index} $($\emph{Gittins index} if the project is nonrestless$).$
\end{definition}

\begin{remark}
\label{re:strongindx} \hspace{1in}
\begin{itemize}
\item[\textup{(i)}] \citet{whit88b} defined indexability when $c(x, a) \triangleq a$   
through an equivalent problem collection $\mathcal{P}' \triangleq \{P_\lambda'\colon \lambda \in \mathbb{R}\}$ differing from $\mathcal{P}$ (see (\ref{eq:lpricepfg1})) in the net rewards, which are $r(x, a) + \lambda (1-a)$ for $P_\lambda'$ (so $\lambda$ represents a \emph{subsidy for passivity}) and $r(x, a) - \lambda a$ for $P_\lambda$ (so $\lambda$ is a \emph{price of activity}). 
Yet, their equivalence is not ensured for projects with general  $c(x, a)$,
where net rewards for $P_\lambda'$ and $P_\lambda$ would be 
$r(x, a) + \lambda c(x, 1-a)$ and $r(x, a) - \lambda c(x, a)$, respectively. 
In such a general case, \citet{nmmp02} introduced a definition of indexability through the problem collection $\mathcal{P}$ as herein. The resulting index is a \emph{generalized Whittle index}, or a \emph{generalized Gittins index} if projects are nonrestless. 
\item[\textup{(ii)}] 
\citet[p.\ 291]{whit88b} defined indexability as a property of the family of $P_\lambda'$-\emph{optimal passive sets} $D(\lambda)$, in that $D(\lambda)$ ``\textit{increases monotonically from $\emptyset$ to $\mathsf{X}$ as $\nu$} (his notation for our $\lambda$) \textit{increases from $-\infty$ to $+\infty$}.'' He defined the index as the root (assumed to be unique) in $\lambda$ of (\ref{eq:bothactopt}). 
\item[\textup{(iii)}] The project is indexable with index $\lambda^*(\cdot)$ if the $P_\lambda$-optimal action in each state $x$ is nonincreasing in $\lambda$ with both actions being $P_\lambda$-optimal only for  $\lambda = \lambda^*(x);$  thus$,$ for $\lambda < \lambda^*(x)$ $($resp.\ $\lambda > \lambda^*(x))$, $a = 1$ $($resp.\ $a = 0)$ is the only $P_\lambda$-optimal action (cf.\ \citet[section 2.3]{nmasmta14}); equivalently, 
 for each $x,$ equation (\ref{eq:bothactopt}) in $\lambda$ has the unique root $\lambda^*(x)$,  and $\Delta_{a=0}^{a=1} V_\lambda^{\langle a, *\rangle}(x) > 0$ for $\lambda < \lambda^*(x)$ and $\Delta_{a=0}^{a=1} V_\lambda^{\langle a, *\rangle}(x) < 0$ for $\lambda > \lambda^*(x).$ 
\end{itemize}
\end{remark}

\subsection{Thresholdability and threshold-indexability: weak and strong versions.}
\label{s:assub}
We address indexability in the case that it holds along with optimality of
threshold policies. 
For any threshold $z$ we define 
the \emph{$z$-policy}, with active and passive regions $S^{1}(z) \triangleq \{x \in \mathsf{X}\colon x > z\}$ and
$S^{0}(z) \triangleq \{x \in \mathsf{X}\colon x \leqslant z\}$; and 
 the  \emph{$z^{\scriptscriptstyle -}$-policy}, with active and passive regions $S^{1}(z^{\scriptscriptstyle -}) \triangleq \{x \in \mathsf{X}\colon x \geqslant z\}$ and
$S^{0}(z^{\scriptscriptstyle -}) \triangleq \{x \in \mathsf{X}\colon x < z\}$. 
  We   write their 
metrics as $F(\nu_0, z)$, $G(\nu_0, z)$, $F(\nu_0, z^{\scriptscriptstyle -})$ and $G(\nu_0, z^{\scriptscriptstyle -})$. 

The notation $F(\nu_0, \cdot)$ and $G(\nu_0,\cdot)$ refers to metrics $F(\nu_0, z)$ and $G(\nu_0, z)$ as functions of the threshold $z$ on $\mathbb{R}$.
The ambiguity of the notation $F(\nu_0, z^{\scriptscriptstyle -})$ and $G(\nu_0, z^{\scriptscriptstyle -})$ is resolved in
 Lemma \ref{pro:cadlagFG}(a), ensuring that the latter metrics are also the left limits at $z$ of $F(\nu_0, \cdot)$ and $G(\nu_0,\cdot)$, respectively. 
 
To formalize the idea of optimality of threshold policies, we distinguish between weak and strong versions of such a concept. 
The following weak version corresponds to the case that $P_\lambda$-optimal actions are partially characterized by a threshold. Note that the term \emph{thresholdability} has been used in a similar sense in the literature. See, e.g.,  \citet{ouyangetal15}.
 
\begin{definition}[Thresholdability]
\label{def:potp}
We call the project \emph{thresholdable}
if  
there is a map $z^*\colon \mathbb{R} \to \overline{\mathbb{R}}$ such that, 
for every project state $x \in \mathsf{X}$ and  price $\lambda \in \mathbb{R}$:
\begin{equation}
\label{eq:thrlsigndef}
\begin{split}
\textup{if } \, x \leqslant z^*(\lambda) & \, \textup{ then } \, 
 \Delta_{a=0}^{a=1} V_\lambda^{\langle a, *\rangle}(x) \leqslant 0; \quad \textup{ and } \\
\textup{if } \, x \geqslant z^*(\lambda)  & \, \textup{ then } \, 
 \Delta_{a=0}^{a=1} V_\lambda^{\langle a, *\rangle}(x) \geqslant 0.
\end{split}
\end{equation}
We call any such $z^*(\cdot)$ an \emph{optimal threshold map}.
\end{definition}

\begin{remark}
\label{re:thrlsigndef} 
\begin{itemize}
\item[\textup{(i)}] If the project is thresholdable with $z^*(\cdot)$ satisfying (\ref{eq:thrlsigndef}), the optimal value function is $V_\lambda^*(x) = V_\lambda(x, z^*(\lambda)) = V_\lambda(x, z^*(\lambda)^{\scriptscriptstyle -})$ for every $x$ and $\lambda$.
\item[\textup{(ii)}]
The project is thresholdable if for every price $\lambda \in \mathbb{R}$ there is a threshold $z \in \overline{\mathbb{R}}$ (not necessarily unique) such that 
both the $z$-policy and the $z^{\scriptscriptstyle -}$-policy are $P_\lambda$-optimal.
\item[\textup{(iii)}] A sufficient condition for thresholdability is that, for every price $\lambda \in \mathbb{R}$,  $V_\lambda^{\langle a, *\rangle}(x)$ is \emph{supermodular} in $(x, a)$, i.e., $\Delta_{a=0}^{a=1} V_\lambda^{\langle a, *\rangle}(x)$ is nondecreasing in $x$.
\item[\textup{(iv)}] Suppose a project is thresholdable. Could we use such information to establish its indexability, using only intrinsic properties of its optimal threshold maps $z^*(\cdot)$? Clearly the answer is negative, since (\ref{eq:thrlsigndef}) does not fully characterize when each action is optimal, whereas (\ref{eq:indxsigndef}) does.
\end{itemize}
\end{remark}

Motivated by Remark \ref{re:thrlsigndef}(iv), we will further consider the following strong version of the thresholdability property, under which the $P_\lambda$-optimal actions are fully characterized by a threshold.

\begin{definition}[Strong thresholdability]
\label{def:spotp}
We call the project \emph{strongly thresholdable}
if  
there is a map $z^*\colon \mathbb{R} \to \overline{\mathbb{R}}$ such that
\begin{equation}
\label{eq:sthrlsigndef} 
\sgn \Delta_{a=0}^{a=1} V_\lambda^{\langle a, *\rangle}(x) = \sgn (x - z^*(\lambda)), \quad \textup{for } x \in \mathsf{X} \textup{  and  } \lambda \in \mathbb{R}.
\end{equation}
\end{definition}

\begin{remark}
\label{re:strongthr} \hspace{1in}
\begin{itemize}
\item[\textup{(i)}] There can only be one optimal threshold map $z^*(\cdot)$ satisfying (\ref{eq:sthrlsigndef}).
\item[\textup{(ii)}] 
The project is strongly thresholdable with optimal threshold map $z^*(\cdot)$ if the $P_\lambda$-optimal action is nondecreasing in the state $x$, with both actions being $P_\lambda$-optimal only for  $x = z^*(\lambda);$  thus$,$ for $x < z^*(\lambda)$ $($resp.\ $x > z^*(\lambda))$, $a = 0$ $($resp.\ $a = 1)$ is the only $P_\lambda$-optimal action (cf.\ Remark \ref{re:strongindx}(iii)); equivalently, 
 for each $\lambda,$ equation (\ref{eq:bothactopt}) in $x$ has $z^*(\lambda)$ as its unique root  and, further, $\Delta_{a=0}^{a=1} V_\lambda^{\langle a, *\rangle}(x) < 0$ for $x < z^*(\lambda)$ and $\Delta_{a=0}^{a=1} V_\lambda^{\langle a, *\rangle}(x) > 0$ for $x > z^*(\lambda);$ or,  
 for every  $\lambda$, the $P_\lambda$-optimal active and passive regions are given by
\begin{equation}
\label{eq:tSstar10}
S^{*, 1}_\lambda = \{x \in \mathsf{X}\colon x \geqslant z^*(\lambda)\} \; \textup{ and } \; 
 S^{*, 0}_\lambda = \{x \in \mathsf{X}\colon x \leqslant z^*(\lambda)\}.
\end{equation}
\end{itemize}
\end{remark}

The following definition formalizes the concept of indexability holding consistently with optimality of threshold policies. 

\begin{definition}[Threshold-indexability]
\label{def:indxtp}
We call the project \emph{threshold-indexable} if it is thresholdable and indexable.
\end{definition}

Our main aim is to identify sufficient conditions for threshold-indexability.
We will further consider the corresponding strong version of the latter property.

\begin{definition}[Strong threshold-indexability]
\label{def:sindxtp}
We call the project \emph{strongly threshold-indexable} if it is strongly thresholdable and indexable.
\end{definition}

\subsection{Marginal metrics and the MP index.}
\label{s:mrrum} As key tools in our indexability analyses, 
we will use the concepts of \emph{marginal metrics} and \emph{marginal productivity} (MP) index,  introduced in \citet{nmmp02} for the analysis of discrete-state projects. 
 
For an action $a$ and policy $\pi$,  let $\langle a, \pi\rangle$ be the policy that takes action
$a$ in period $t = 0$ and adopts 
$\pi$ afterwards.
We measure the  changes in reward earned and in resource expended, respectively, resulting from a change from $\langle 0, \pi\rangle$  to $\langle 1, \pi\rangle$ starting from $x$, by the
\emph{marginal reward metric} 
\begin{equation}
\label{eq:mrm1}
f(x, \pi)\triangleq  \diff\limits_{a=0}^{a=1}  F(x, \langle a, \pi\rangle) = F(x, \langle 1, \pi\rangle) - F(x, \langle 0, \pi\rangle)
\end{equation}
and  the \emph{marginal resource $($\!usage$)$ metric} 
\begin{equation}
\label{eq:mwm1}
g(x, \pi)\triangleq \diff\limits_{a=0}^{a=1} G(x, \langle a, \pi\rangle).
\end{equation}

\begin{remark}
\label{re:mpmprop}
\begin{itemize}
\item[\textup{(i)}] It follows  from Remark \textup{\ref{re:assboundfs}(ii)}  and \textup{(\ref{eq:mrm1}, \ref{eq:mwm1})} that$,$ for any policy $\pi \in \Pi,$  
\begin{equation}
\label{eq:fgxpibnd}
f(\cdot, \pi), g(\cdot, \pi) \in \mathbb{B}_w(\mathsf{X}), \textup{ with } \max \{\| f(\cdot, \pi)\|_w , \| g(\cdot, \pi)\|_w\} \leqslant   2 M_\gamma.
\end{equation}
\item[\textup{(ii)}] For a given active region
$B,$ the functions $f(\cdot, B)$ and $g(\cdot, B)$ can be evaluated by
\begin{equation*}
f(x, B) = 
\begin{cases}
\displaystyle F(x, B) - r(x, 0) - \beta \int F(y, B) \,
  \kappa^{0}(x, dy), &
  \quad x \in B \\
\displaystyle r(x, 1) + \beta \int F(y, B) \,
  \kappa^{1}(x, dy) - F(x, B), & \quad x \in B^c,
\end{cases}
\end{equation*}
\begin{equation*}
g(x, B) = 
\begin{cases}
\displaystyle G(x, B) - c(x, 0) - \beta \int G(y, B) \,
  \kappa^{0}(x, dy), &
  \quad 
x \in B \\
\displaystyle c(x, 1) + \beta \int G(y, B) \,
  \kappa^{1}(x, dy) - G(x, B), & \quad x \in B^c;
\end{cases}
\end{equation*}
or$,$ in terms of the linear operator $\mathcal{L}^\star$ in $(\ref{eq:tbetastu}),$
\begin{equation}
\label{eq:fxBeq}
f(x, B) = 
\begin{cases}
\displaystyle  \mathcal{L}^\star F(\cdot, B)(x, 0) - r(x, 0), &
  \quad 
x \in B \\
\displaystyle r(x, 1) - \mathcal{L}^\star F(\cdot,
  B)(x, 1), & \quad x \in B^c,
\end{cases}
\end{equation}
\begin{equation}
\label{eq:gxBeq}
g(x, B) = 
\begin{cases}
\displaystyle  \mathcal{L}^\star G(\cdot, B)(x, 0) - c(x, 0), &
  \quad 
x \in B \\
\displaystyle c(x, 1) - \mathcal{L}^\star G(\cdot,
  B)(x, 1), & \quad x \in B^c.
\end{cases}
\end{equation}
\item[\textup{(iii)}]  For a given project state $x$ and policy
$\pi,$ $f(x, \pi)$ and $g(x, \pi)$ can be represented as
\begin{equation}
\label{eq:fxBgxB}
f(x, \pi) = \diff\limits_{a=0}^{a=1} \{r(x, a) + \beta F(\kappa^a(x, \cdot), \pi)\} \quad \textup{and} \quad
g(x, \pi) = \diff\limits_{a=0}^{a=1} \{c(x, a) + \beta G(\kappa^a(x, \cdot), \pi)\}.
\end{equation}
\end{itemize}
\end{remark}

If $g(x, B)\neq 0$ for a given $x$ and $B$, we measure the MP of engaging the project in state $x$ at time $t = 0$, given that the $B$-policy is adopted thereafter, by the
\emph{MP metric} 
\begin{equation}
\label{eq:mpm1}
m(x, B)\triangleq \frac{f(x, B)}{g(x, B)}.
\end{equation}

We  shall further write as
$f(x, z)$, $g(x, z)$, $m(x, z)$ and $f(x, z^{\scriptscriptstyle -})$, $g(x, z^{\scriptscriptstyle -})$, $m(x, z^{\scriptscriptstyle -})$ the  marginal metrics for the $z$-policy and the $z^{\scriptscriptstyle -}$-policy, respectively (see Section \ref{s:assub}).

We next use the  MP metric to define the  \emph{MP index}.

\begin{definition}[MP index]
\label{def:mpi}
Suppose $g(x, x) \neq 0$ for any state $x.$ The project's \emph{MP index} is the map $m\colon \mathsf{X} \to \mathbb{R}$ given by 
\begin{equation}
\label{eq:phimpx}
m(x) \triangleq m(x, x) = \frac{f(x, x)}{g(x, x)}, \quad x \in \mathsf{X}.
\end{equation} 
\end{definition}

Note the abuse of notation in (\ref{eq:phimpx}), since $m$ is used both for  the MP index and the MP metric. 

\subsection{Main result: a verification theorem for threshold-indexability.}
\label{s:pclbsic}
We next present our main result, giving  sufficient conditions for threshold-indexability.
The conditions correspond to the concept of \emph{PCL-indexability}, which is extended here 
 to the real-state setting.

\begin{definition}[PCL-indexability]
\label{def:pcli}  
\rm
We call the project \emph{PCL-indexable}  (with respect to threshold policies) if the
following \emph{PCL-indexability conditions} hold:
\begin{itemize}
\item[\textup{(PCLI1)}] the marginal resource metric satisfies 
$g(x, z)  > 0$ for every $x \in \mathsf{X}$ and $z \in \overline{\mathbb{R}};$ 
\item[\textup{(PCLI2)}] the MP index $m(\cdot)$ is nondecreasing and  continuous on $\mathsf{X};$ 
\item[\textup{(PCLI3)}] for each state $x,$ the metrics $F(x, \cdot),$ $G(x, \cdot)$ and $m(\cdot)$ are related by (\ref{eq:pcli3}).
\end{itemize}
\end{definition}

\begin{remark}
\label{re:pcliii}
\begin{itemize}
\item[\textup{(i)}] Conditions \textup{(PCLI1, PCLI2)} match those given for discrete-state projects in \textup{\citet{nmaap01,nmmp02,nmmor06}}$,$ except for the continuity in (PCLI2)$.$
\item[\textup{(ii)}] Under (PCLI1), the MP metric $m(x, z)$ (see (\ref{eq:mpm1})) is well defined for every $x$ and $z.$ 
\item[\textup{(iii)}] Under (PCLI1, PCLI2), we extend the MP index domain to  $\overline{\mathbb{R}}$ by 
$m(z) \triangleq m(\ell)$ for $z < \ell$ if $\ell > -\infty$, $m(z) \triangleq m(u)$ for $z > u$ if $u < \infty$, 
$m(-\infty) \triangleq \lim_{z \to -\infty} m(z)$, 
and $m(\infty) \triangleq \lim_{x \to \infty} m(x)$.
\item[\textup{(iv)}] Under (PCLI1, PCLI2), it follows immediately that the range $m(\mathsf{X})$ of $m(\cdot)$ is an interval.
\item[\textup{(v)}] (PCLI3) requires that, for each state $x$, $F(x, \cdot)$ be an  indefinite LS integral of the MP index $m(\cdot)$ with respect to $G(x, \cdot).$ 
The existence of such integrals is justified in Lemma \ref{lma:eaafpcli3}.
\item[\textup{(vi)}] For discrete-state projects$,$  (PCLI1) implies \textup{(PCLI3)}. See \citet[theorem 6.4(b)]{nmmp02} for the finite-state case and \citet[lemma 5.5.(a)]{nmmor06} for the denumerable-state case. 
\item[\textup{(vii)}] For any initial-state distribution $\nu_0 \in \mathbb{P}_w(\mathsf{X}),$ \textup{(PCLI3)} yields (cf.\ Remark \ref{re:assboundfs}(iii))
\begin{equation}
\label{eq:pcli3nu}
F(\nu_0, z_2) -  F(\nu_0, z_1) = \int_{(z_1, z_2]} m(z) \, G(\nu_0, dz), \quad -\infty < z_1 < z_2   < \infty.
\end{equation}
\end{itemize}
\end{remark}

For a project satisfying (PCLI1, PCLI2), define for any price $\lambda \in \mathbb{R}$ the  \emph{threshold set}
\begin{equation}
\label{eq:tslambm}
\mathsf{Z}_\lambda \triangleq 
\begin{cases}
\{z \in \overline{\mathbb{R}}\colon m(z) = \lambda\} & \textup{ if } \, \{x \in \mathsf{X}\colon m(x) = \lambda\} \neq \emptyset \\
\{z \in \overline{\mathbb{R}}\colon z < \ell\} & \textup{ if } \, \ell > -\infty \, \textup{ and } \, \lambda < m(x) \, \textup{ for every } \, x \in \mathsf{X} \\
\{-\infty\} & \textup{ if } \, \ell = -\infty \, \textup{ and } \, \lambda < m(x) \, \textup{ for every } \, x \in \mathsf{X} \\
\{z \in \overline{\mathbb{R}}\colon z > u\} & \textup{ if } \, u < \infty \, \textup{ and } \, \lambda > m(x) \, \textup{ for every } \, x \in \mathsf{X} \\
\{\infty\} & \textup{ if } \, u = \infty \, \textup{ and } \, \lambda > m(x) \, \textup{ for every } \, x \in \mathsf{X}.
\end{cases}
\end{equation}

Under (PCLI1, PCLI2), we shall refer to the 
\emph{generalized inverse} (cf.\ \citet[p.\ 413]{falkTeschl12}) of the MP index.
We shall call a threshold map $\zeta\colon \mathbb{R} \to \overline{\mathbb{R}}$ a \emph{generalized inverse} of $m(\cdot)$ if 
\begin{equation}
\label{eq:gimdef}
\zeta(\lambda) \in \mathsf{Z}_\lambda, \quad \lambda \in \mathbb{R}.
\end{equation}

The following result gives two key properties of such generalized inverses of the MP index.

\begin{lemma}
\label{lma:geninvprop}
Assume \textup{(PCLI1, PCLI2)} hold and 
let $\zeta(\cdot)$ be a generalized inverse of $m(\cdot).$  
For any project state $x$ and resource price $\lambda$ the following holds\textup{:}
\begin{itemize}
\item[\textup{(a)}] if $x \leqslant \zeta(\lambda)$ then $m(x) \leqslant \lambda;$ and if $x \geqslant \zeta(\lambda)$ then $m(x) \geqslant \lambda;$
\item[\textup{(b)}] $\zeta(\cdot)$ is increasing on the interval $m(\mathsf{X}).$
\end{itemize}
\end{lemma}
\proof
(a) This part follows follows straightforwardly  from (\ref{eq:tslambm}) and (\ref{eq:gimdef}).

(b) Let $\zeta(\cdot)$ be a generalized inverse of  $m(\cdot)$ and 
suppose $\zeta(\cdot)$ were not increasing on $m(\mathsf{X}),$ which is an interval by (PCLI2). Then there would exist $m(x_1) = \lambda_1 < \lambda_2 = m(x_2)$ with $x_1, x_2 \in \mathsf{X}$ and $\zeta(\lambda_1) \geqslant \zeta(\lambda_2)$. Now, from  $\zeta(\lambda_k) \in \mathsf{Z}_{\lambda_k}$, (\ref{eq:tslambm}), (\ref{eq:gimdef}) and nondecreasingness of $m(\cdot)$ it would follow that $\lambda_1 = m(\zeta(\lambda_1)) \geqslant m(\zeta(\lambda_2)) = \lambda_2$, a contradiction. Hence $\zeta(\cdot)$ must be increasing on $m(\mathsf{X}).$ 
 \endproof

We next state the verification theorem.

\begin{theorem}
\label{the:pcliii}
Under $($\textup{PCLI1}--\textup{PCLI3}$),$ 
the project is threshold-indexable with Whittle index 
$m(\cdot),$ and its optimal threshold maps are  the generalized inverses $\zeta(\cdot)$ of $m(\cdot).$
\end{theorem}

The proof of Theorem \ref{the:pcliii}, given in Section \ref{s:ptpcliii}, builds on a number of implications of PCL-indexability conditions (PCLI1--PCLI3), which are developed in Section \ref{s:saom}--\ref{s:opicfpmpm}.

\section{Motivation of the PCL-indexability conditions.}
\label{s:potip}
This section presents results motivating the PCL-indexability conditions.

We start by characterizing thresholdability for an indexable project, in its weak and strong versions (see Definitions \ref{def:potp} and \ref{def:spotp}),  in terms  of properties of the Whittle index. 
The following result shows that condition (PCLI2) is necessary for threshold-indexability, and further characterizes optimal threshold maps as generalized inverses of the Whittle index (cf.\ Theorem \ref{the:pcliii}).

\begin{lemma}
\label{lma:iwrttpndlc}
Suppose the project is indexable with Whittle index $\lambda^*(\cdot).$ Then
it is thresholdable  iff $\lambda^*(\cdot)$  is nondecreasing and continuous$,$ in which case its optimal threshold maps $z^*(\cdot)$ are  the generalized inverses of $\lambda^*(\cdot),$ which are increasing on 
$\lambda^*(\mathsf{X}).$
\end{lemma}
\proof 
Suppose the project is thresholdable. To obtain a contradiction, 
suppose $\lambda^*(\cdot)$ were not nondecreasing.  Then there would be states $x < y$ with $\lambda^*(y) < \lambda^*(x)$ and, for any price $\lambda^*(y) < \lambda < \lambda^*(x)$, indexability would imply (see Remark \ref{re:strongindx}(iii)) that  $a = 0$ and $a = 1$ are the only $P_\lambda$-optimal actions in $y$ and  $x$, respectively, which would rule out $P_\lambda$-optimality of a threshold policy. Hence $\lambda^*(\cdot)$ must be nondecreasing. It further follows that $\lambda^*(\cdot)$ can only have jump discontinuities, so the left and right limits $\lambda^*(x^{\scriptscriptstyle -})$ and $\lambda^*(x^{\scriptscriptstyle +})$ exist at states $x > \ell$ and $x < u$, respectively. 

To show that $\lambda^*(\cdot)$ is left-continuous, suppose 
$\lambda^*(y^{\scriptscriptstyle -}) < \lambda^*(y)$ for some $y > \ell$.
Then for any $\lambda^*(y^{\scriptscriptstyle -}) < \lambda < \lambda^*(y)$ and  $x_1 < y < x_2$, it would be $\lambda^*(x_1) \leqslant \lambda^*(y^{\scriptscriptstyle -}) < \lambda < \lambda^*(y) < \lambda^*(x_2)$ by nondecreasingness of $\lambda^*(\cdot)$.
Hence, indexability (cf.\ Remark \ref{re:strongindx}(iii)) would imply that the only $P_\lambda$-optimal threshold policy would be the $y^{\scriptscriptstyle -}$-policy, contradicting thresholdability (cf.\ Remark \ref{re:thrlsigndef}(i)). 

As for right-continuity of $\lambda^*(\cdot)$, if there were a state $y < u$ with $\lambda^*(y) < \lambda^*(y^{\scriptscriptstyle +})$, then
 for any $\lambda^*(y) < \lambda < \lambda^*(y^{\scriptscriptstyle +})$ and states $x_1 < y < x_2$ it would be $\lambda^*(x_1) \leqslant \lambda^*(y) < \lambda < \lambda^*(y^{\scriptscriptstyle +}) \leqslant \lambda^*(x_2)$. Hence by indexability the $y$-policy would be the only $P_\lambda$-optimal threshold policy, again a contradiction. 

As for the reverse implication, suppose  $\lambda^*(\cdot)$ is nondecreasing and continuous, and fix
$\lambda \in \mathbb{R}$. Then, it is immediate from the assumptions that for the following choices of $z$ both the $z$-policy and the $z^{\scriptscriptstyle -}$-policy are $P_\lambda$-optimal:
if $\lambda^*(x) < \lambda$ for every state $x$, take any $z > u$ if $u < \infty$, and $z = \infty$ if $u = \infty$;  if $\lambda^*(x) > \lambda$ for every state $x$, take any $z < \ell$ if $\ell > -\infty$, and $z = -\infty$ if $\ell = -\infty$; and, otherwise, take as $z$ any state with $\lambda^*(z) = \lambda$, whose existence is ensured by (PCLI2).  This shows both that the project is thresholdable (see Remark \ref{re:thrlsigndef}(i)), and that its optimal threshold maps $z^*(\cdot)$ are the generalized inverses of $\lambda^*(\cdot)$. It further follows as in the proof of Lemma \ref{lma:geninvprop}(b) that any such $z^*(\cdot)$ must be  increasing on 
$\lambda^*(\mathsf{X}).$
 \endproof

The next result characterizes via marginal metrics the optimal regions $S^{*, a}_\lambda$ in (\ref{eq:Sstarlam}).
Let $\lambda \in \mathbb{R}$.

\begin{lemma}
\label{lma:mpm} Let policy $\pi^*$ be $P_\lambda$-optimal$.$  Then
\begin{equation*}
\begin{split}
S^{*, 1}_\lambda & = \{x \in \mathsf{X}\colon f(x, \pi^*) - \lambda g(x, \pi^*) \geqslant 0\} \quad \textup{and}  \quad
S^{*, 0}_\lambda = \{x \in \mathsf{X}\colon f(x, \pi^*) - \lambda g(x, \pi^*) \leqslant 0\}.
\end{split}
\end{equation*}
\end{lemma}
\proof
Since $\pi^*$ is $P_\lambda$-optimal, we can reformulate  (\ref{eq:Sstarlam}) as  
\begin{equation}
\label{eq:refBEst}
\begin{split}
S^{*, a}_\lambda  & = \big\{x \in \mathsf{X}\colon V_\lambda(x, \langle a, \pi^*\rangle) \geqslant V_\lambda(x, \langle 1-a, \pi^*\rangle)\big\}, \quad a \in \{0, 1\}.
\end{split}
\end{equation}

Now,  it immediately follows that 
\begin{equation}
\label{eq:refBEpist}
\diff\limits_{a=0}^{a=1} V_\lambda(x, \langle a, \pi^*\rangle) = f(x, \pi^*) - \lambda g(x, \pi^*),
\end{equation}
which, along with (\ref{eq:refBEst}), yields the result. \qquad
 \endproof

The next result justifies using the MP index $m(\cdot)$ as a candidate for 
the project's Whittle index.

\begin{proposition}
\label{lma:phiwphimpi}
Let the project be threshold-indexable with Whittle index $\lambda^*(\cdot).$ Then
\begin{itemize}
\item[\textup{(a)}] $f(x, x) - \lambda^*(x) g(x, x) = 0 = f(x, x^{\scriptscriptstyle -}) - \lambda^*(x) g(x, x^{\scriptscriptstyle -})$ for every state $x;$
\item[\textup{(b)}]
 $\lambda^*(x) = m(x)$ for any state $x$ with $g(x, x) \neq 0;$
\item[\textup{(c)}] under \textup{(PCLI1),}  $\lambda^*(\cdot) = m(\cdot).$
\end{itemize} 
\end{proposition}
\proof
(a) 
Since the project is threshold-indexable,   $\lambda^*(\cdot)$ is nondecreasing by Lemma \ref{lma:iwrttpndlc}. It follows from this and Definition \ref{def:indx} that both the $x$-policy and the $x^{\scriptscriptstyle -}$-policy are $P_\lambda$-optimal for $\lambda = \lambda^*(x)$.  Hence, Lemma \ref{lma:mpm} implies that $f(x, x) - \lambda^*(x) g(x, x) = f(x, x^{\scriptscriptstyle -}) - \lambda^*(x) g(x, x^{\scriptscriptstyle -}) = 0$.

(b) This part follows from the first identity in (a) and (\ref{eq:phimpx}). 

(c) This part follows from (b) and (PCLI1). 
\qquad
 \endproof

The following two lemmas will be used to establish in Proposition \ref{lma:intgradz} the necessity of (PCLI3) for threshold-indexability.
Note that, for any $x$, the map $\lambda \mapsto V_\lambda^*(x)$ is convex on $\mathbb{R}$, as $V_\lambda^*(x)$ is the supremum of linear maps in $\lambda$. Hence, it is differentiable almost everywhere and has a nonempty set 
 of \emph{subgradients} (see, e.g.,  \citet[\S1.5]{NicuPer06}) at every $\lambda$, which we denote by $\partial  V_\lambda^*(x)$. 
The following result shows how to obtain subgradients of such a map.
Let $x \in \mathsf{X}$ and $\lambda \in \mathbb{R}$.

\begin{lemma}
\label{lma:sbgrad} 
If $\pi^*$ is a $P_{\lambda}$-optimal policy$,$ then 
$-G(x, \pi^*) \in \partial  V_\lambda^*(x).$
\end{lemma}
\proof
Let $\delta > 0$. Using that $V_{\lambda}^*(x) = F(x, \pi^*) - \lambda G(x, \pi^*)$, we have
\[
V_{\lambda}^*(x) - \delta G(x, \pi^*) = F(x, \pi^*) - (\lambda + \delta) G(x, \pi^*) = V_{\lambda + \delta}(x, \pi^*) \leqslant V_{\lambda + \delta}^*(x)
\]
and
\[
V_{\lambda}^*(x) + \delta G(x, \pi^*) = F(x, \pi^*) - (\lambda - \delta) G(x, \pi^*) = V_{\lambda - \delta}(x, \pi^*) \leqslant V_{\lambda - \delta}^*(x),
\]
whence 
\[
V_{\lambda}^*(x) - V_{\lambda - \delta}^*(x) \leqslant - \delta G(x, \pi^*) \leqslant V_{\lambda + \delta}^*(x) - V_{\lambda}^*(x).
\]
This shows that $-G(x, \pi^*)$ is a subgradient of $V_{\lambda}^*(x)$.
 \endproof

\begin{lemma}
\label{lma:intgradl} 
Suppose the project is thresholdable and let $z^*(\cdot)$ be an optimal threshold map$.$
Then for any state $x$ and prices $\lambda_1 < \lambda_2,$ 
\begin{itemize}
\item[\textup{(a)}]
$V_{\lambda_1}^*(x) - V_{\lambda_2}^*(x) = \int_{(\lambda_1, \lambda_2]} G(x, z^*(\lambda)) \, d\lambda;$
\item[\textup{(b)}]
$\Delta_{a=0}^{a=1} V_{\lambda_1}^{\langle a, *\rangle}(x) - \Delta_{a=0}^{a=1} V_{\lambda_2}^{\langle a, *\rangle}(x) = \int_{(\lambda_1, \lambda_2]} g(x, z^*(\lambda)) \, d\lambda.$
\end{itemize}
\end{lemma}
\proof (a)
This follows from the integral representation of convex functions (see
\citet[proposition 1.6.1]{NicuPer06}) and Lemma \ref{lma:sbgrad}, as $-G(x, z^*(\lambda)) \in \partial  V_\lambda^*(x)$.

(b) We have
\begin{align*}
V_{\lambda_1}^{\langle a, *\rangle}(x) - V_{\lambda_2}^{\langle a, *\rangle}(x) & = 
(\lambda_2-\lambda_1) c(x, a) + \beta \int \big[V_{\lambda_1}^*(y) - V_{\lambda_2}^*(y)\big] \, \kappa^a(x, dy) \\
& = (\lambda_2-\lambda_1) c(x, a) + \beta \int \int_{(\lambda_1, \lambda_2]} G(y, z^*(\lambda)) \, d\lambda \, \kappa^a(x, dy) \\
& = \int_{(\lambda_1, \lambda_2]} \bigg[c(x, a) + \beta \int G(y, z^*(\lambda)) \, \kappa^a(x, dy)\bigg] \,  d\lambda,
\end{align*}
where we have used  that $V_{\lambda}^{\langle a, *\rangle}(x) = r(x, a) - \lambda c(x, a) + \beta \int V_{\lambda}^*(y) \, \kappa^a(x, dy)$, part (a), and Fubini's theorem (using that $G(x, \pi) \geqslant 0$).
Therefore we obtain, as required,
\[
\Delta_{a=0}^{a=1} \big[V_{\lambda_1}^{\langle a, *\rangle}(x) - V_{\lambda_2}^{\langle a, *\rangle}(x)\big] = \int_{(\lambda_1, \lambda_2]} g(x, z^*(\lambda)) \,  d\lambda.
\]
\endproof

The following result shows that condition (PCLI3) in (\ref{eq:pcli3}) is necessary for threshold-indexability.
Its proof relies on $G(x, \cdot)$ being a right-continuous function of bounded variation on finite intervals, to ensure the validity of identities (\ref{eq:chvar1ft12}) and (\ref{eq:intparts1}).
This is proved in Lemmas \ref{pro:cadlagFG} and \ref{pro:Gnondecr}.

\begin{proposition}
\label{lma:intgradz} 
Suppose the project is threshold-indexable with index $\lambda^*(\cdot).$ Then for any $x,$
\[
F(x, z_2) - F(x, z_1) = \int_{(z_1, z_2]} \lambda^*(z) \, G(x, dz), \quad -\infty < z_1 < z_2   < \infty.
\]
\end{proposition}
\proof Assume with no loss of generality that $z_1, z_2 \in \mathsf{X}$.
Let $z^*(\cdot)$ be an optimal threshold map. By Lemma \ref{lma:iwrttpndlc}, $z^*(\cdot)$ is a generalized inverse of the nondecreasing continuous map $\lambda^*(\cdot)$, and hence we can apply the change of variable formula in 
 \citet[Eq.\ (1)]{falkTeschl12} to obtain
\begin{equation}
\label{eq:chvar1ft12}
\int_{(z_1, z_2]} G(x, z) \, \lambda^*(dz) = \int_{(\lambda_1, \lambda_2]} G(x, z^*(\lambda)) \, d\lambda,
\end{equation}
where $\lambda_k = \lambda^*(z_k)$ for $k = 1, 2$.

On the other hand, integration by parts (see \citet[theorem 6.2.2]{cartvanBrunt00}) gives
\begin{equation}
\label{eq:intparts1}
\int_{(z_1, z_2]} G(x, z) \, \lambda^*(dz) = \lambda_2 G(x, z_2)  - \lambda_1 G(x, z_1)  - \int_{(z_1, z_2]} \lambda^*(z) \, G(x, dz). 
\end{equation}

Furthermore, we can write 
\begin{equation}
\label{eq:vlambda12x}
V_{\lambda_1}^*(x) - V_{\lambda_2}^*(x) = 
F(x, z_1) - \lambda_1 G(x, z_1) - F(x, z_2) + \lambda_2 G(x, z_2).
\end{equation}

Now, we use in turn (\ref{eq:vlambda12x}), Lemma \ref{lma:intgradl}(a), (\ref{eq:intparts1}) and (\ref{eq:chvar1ft12}) to obtain
\begin{align*}
F(x, z_2) - F(x, z_1) & = \lambda_2 G(x, z_2) - \lambda_1 G(x, z_1)  - V_{\lambda_1}^*(x) + V_{\lambda_2}^*(x) \\
& = \int_{(z_1, z_2]} G(x, z) \, \lambda^*(dz) +  \int_{(z_1, z_2]} \lambda^*(z) \, G(x, dz) - \int_{(\lambda_1, \lambda_2]} G(x, z^*(\lambda)) \, d\lambda \\
& = \int_{(z_1, z_2]} \lambda^*(z) \, G(x, dz).
\end{align*}
 \endproof

\section{Indexability of strongly thresholdable projects.}
\label{s:srsti}
We next present results relevant to prevailing approaches to indexability (see Section \ref{s:rmabprsp}), which rely on first
establishing strong thresholdability.

\begin{lemma}
\label{lma:siwrttpndlc}
An indexable project is strongly thresholdable  iff its index $\lambda^*(\cdot)$  is increasing and continuous$,$
in which case its optimal threshold map $z^*(\cdot)$ is increasing and continuous on $\lambda^*(\mathsf{X}).$ 
\end{lemma}
\proof 
Suppose the project is strongly thresholdable. Then, by Lemma \ref{lma:iwrttpndlc}, $\lambda^*(\cdot)$ is nondecreasing and continuous. To show that its Whittle index $\lambda^*(\cdot)$ must be increasing note that,
 from Definitions \ref{def:indx} and \ref{def:spotp}, its unique (see Remark \ref{re:strongthr}(i)) optimal threshold map $z^*(\cdot)$ and Whittle index satisfy \begin{equation}
\label{eq:lmathrlsigndef}
\sgn (x - z^*(\lambda)) = \sgn (\lambda^*(x) - \lambda), \quad \textup{for } x \in \mathsf{X} \textup{  and  } \lambda \in \mathbb{R}.
\end{equation}
Take states $x_1 < x_2$.  Then, writing $\lambda_k \triangleq \lambda^*(x_k)$, (\ref{eq:lmathrlsigndef}) yields $z^*(\lambda_k) = x_k$ for $k = 1, 2$.
We now use again (\ref{eq:lmathrlsigndef}) to obtain, as required,
\[
\sgn (\lambda_1 - \lambda_2) = \sgn (\lambda^*(x_1) - \lambda_2) = \sgn (x_1 - z^*(\lambda_2)) = \sgn (x_1 - x_2) < 0.
\]

Consider now the reverse implication, so $\lambda^*(\cdot)$  is assumed increasing and continuous.
Then it has a unique inverse $z^{*}(\cdot) \triangleq \lambda^{*, -1}(\cdot)$ on the interval $\lambda^*(\mathsf{X})$, which is also increasing and continuous, and satisfies  the identity in (\ref{eq:lmathrlsigndef}) for $x \in \mathsf{X}$ and $\lambda \in \lambda^*(\mathsf{X})$. If $\lambda^*(\mathsf{X}) \neq \mathbb{R}$, it is immediate to extend the domain of such a $z^{*}(\cdot)$ to $\mathbb{R}$ consistently with (\ref{eq:lmathrlsigndef}). This yields strong thresholdability.
 \endproof

The next result shows that to establish indexability of a strongly thresholdable project it does not suffice to prove that its optimal threshold map is increasing. It needs also to be proved that the optimal threshold map is continuous. Note that $z^{*, -1}(\mathsf{X}) \triangleq \{\lambda\colon z^*(\lambda) \in \mathsf{X}\}$.

\begin{proposition}
\label{pro:iwrttpndlc}
Assume the project is strongly thresholdable with optimal threshold map $z^*(\cdot).$
Then it is indexable  iff $z^*(\cdot)$ is increasing and continuous on $z^{*, -1}(\mathsf{X}),$ in which case its Whittle index is also increasing and continuous$,$ being $\lambda^*(x) = z^{*, -1}(x)$ for $x \in \mathsf{X}.$
\end{proposition}
\proof
If the project is indexable, the stated properties of $z^*(\cdot)$ follow by 
Lemma \ref{lma:iwrttpndlc}.

As for the reverse implication, suppose the project is strongly thresholdable with $z^*(\cdot)$ increasing and continuous on $z^{*, -1}(\mathsf{X})$.
Then the restriction of $z^*(\cdot)$ on the domain $z^{*, -1}(\mathsf{X})$ has a unique inverse $\lambda^*(\cdot) \triangleq z^{*, -1}(\cdot)$, which is also increasing and continuous, and satisfies (\ref{eq:lmathrlsigndef}).
This implies that the project is indexable with Whittle index $\lambda^*(\cdot)$, which completes the proof.
 \endproof

\section{LP reformulation of the $\lambda$-price problem starting from $X_0 \sim \nu_0$.}
\label{s:saom}
This section lays down further groundwork for proving Theorem \ref{the:pcliii}.
It draws upon standard results to reformulate the problem 
\begin{equation}
\label{eq:lambdapp}
P_{\lambda}(\nu_0)\colon \qquad \maxim_{\pi \in \Pi} \, V_\lambda(\nu_0, \pi),
\end{equation}
which is the $\lambda$-price problem $P_\lambda$ in (\ref{eq:lpricepfg1}) starting from $X_0 \sim \nu_0 \in \mathbb{P}_w(\mathsf{X})$, as an LP problem. 

For any policy $\pi \in \Pi$, consider the  \emph{state-action occupation measure} $\mu_{\nu_0}^{\pi}$ defined by
\begin{equation}
\label{eq:saom}
\mu_{\nu_0}^{\pi}(\Gamma) \triangleq \Ex_{\nu_0}^\pi\Bigg[\sum_{t=0}^\infty 
 \beta^t 1_{\Gamma}(X_t, A_t)\Bigg] = \sum_{t=0}^\infty \beta^t
\Prob_{\nu_0}^\pi\{(X_t, A_t) \in \Gamma\}, \quad \Gamma \in \mathcal{B}(\mathsf{K}).
\end{equation}
Since the objective of problem $P_{\lambda}(\nu_0)$ can be represented in terms of $\mu_{\nu_0}^{\pi}$ as
\[
V_\lambda(\nu_0, \pi) = \int (r - \lambda c) \,  d\mu_{\nu_0}^{\pi} = 
\int \{r(y, a) - \lambda c(y, a)\} \,  \mu_{\nu_0}^{\pi}(d(y, a)),
\]
we can reformulate the dynamic optimization problem $P_{\lambda}(\nu_0)$ as the static optimization problem
\begin{equation}
\label{eq:psop}
L_{\lambda}(\nu_0)\colon \qquad \maxim_{\mu \in \mathcal{M}_{\nu_0}} \,\int (r - \lambda c) \, d\mu,
\end{equation}
where $\mathcal{M}_{\nu_0} \triangleq \{\mu_{\nu_0}^{\pi}\colon \pi \in \Pi\}$, which is 
 characterized  by linear constraints (see, e.g.,  \citet{heilmann79} and \citet[\S6.3]{herlerLass96}) as outlined below.

The $\mu_{\nu_0}^{\pi}$ belong (see \citet[proposition 7.2.2]{herlerLass99}) to the Banach space  $\mathbb{M}_w(\mathsf{K})$  of
finite \emph{signed measures} $\mu$ on $\mathcal{B}(\mathsf{K})$ with finite $w$-norm (with $|\mu|$ the \emph{total variation} of $\mu$)
\begin{equation}
\label{eq:mubnorm}
\| \mu \|_w \triangleq \sup_{h \in \mathbb{B}_w(\mathsf{K}), \| h \|_w \leqslant 1} \bigg|\int  h \, d\mu\bigg| = \int  w \, d|\mu| = \int  w(y) \, |\mu|(d(y, a)).
\end{equation}
  
For any $\mu \in \mathbb{M}_w(\mathsf{K})$ consider the marginal signed measures on $\mathcal{B}(\mathsf{X})$ defined by 
$\tilde{\mu}^a(S) \triangleq \mu(S \times \{a\})$ for $a \in \{0, 1\}$, and 
$\tilde{\mu}(S) \triangleq \mu(S \times \{0, 1\}) = \sum_{a \in \{0,
  1\}} \tilde{\mu}^a(S)$, which belong to
the Banach space $\mathbb{M}_w(\mathsf{X})$    of finite signed
measures $\nu$ on $\mathcal{B}(\mathsf{X})$ with finite $w$-norm (cf.\ (\ref{eq:pmwnorm}))
\begin{equation}
\label{eq:nubnorm}
\| \nu \|_w \triangleq \sup_{v \in \mathbb{B}_w(\mathsf{X}), \| v \|_w \leqslant 1} \bigg|\int  v \, d\nu\bigg| =
\int  w \, d|\nu| = \int  w(y) \, |\nu|(dy).
\end{equation}
In particular, the \emph{marginal action/state-occupation measures} for $\mu_{\nu_0}^\pi$,  $\tilde{\mu}_{\nu_0}^{\pi, a}$ and $\tilde{\mu}_{\nu_0}^{\pi}$, are in $\mathbb{M}_w(\mathsf{X})$.
 
Both  $(\mathbb{M}_w(\mathsf{K}),
\mathbb{B}_w(\mathsf{K}))$  and  $(\mathbb{M}_w(\mathsf{X}), \mathbb{B}_w(\mathsf{X}))$ are
known to be \emph{dual pairs of vector spaces} (see \citet[\S12.2.A]{herlerLass99}) with respect to the bilinear forms
\begin{align}
\langle \mu, h\rangle \triangleq \int h \, d\mu & = \int h(y, a) \, \mu(d(y, a)), \quad \mu \in \mathbb{M}_w(\mathsf{K}), h \in
\mathbb{B}_w(\mathsf{K}), \label{eq:bfmuu} \\
\langle \nu, v\rangle \triangleq  \int  v \, d\nu & = \int
v(y) \, \nu(dy), \quad \nu \in \mathbb{M}_w(\mathsf{X}), v \in \mathbb{B}_w(\mathsf{X}). \label{eq:bfnuv}
\end{align}

Consider now the linear operators $\mathcal{L}^a\colon \mathbb{M}_w(\mathsf{K}) \to \mathbb{M}_w(\mathsf{X})$ for $a \in \{0, 1\}$, mapping
$\mu$ to 
\begin{equation}
\label{eq:tbadef}
\mathcal{L}^a \mu(S) \triangleq \tilde{\mu}^a(S) - \beta \int
\kappa^a(y, S) \, \tilde{\mu}^a(dy), \quad S \in \mathcal{B}(\mathsf{X}),
\end{equation}
and the linear
operator $\mathcal{L} \triangleq \mathcal{L}^0 + \mathcal{L}^1$. Note that
\begin{equation}
\label{eq:loL}
\mathcal{L} \mu(S) \triangleq \tilde{\mu}(S) - \beta 
\int \kappa^a(y, S) \, \mu(d(y, a)), \quad 
S \in \mathcal{B}\big(\mathsf{X}\big).
\end{equation}
Such operators are \emph{bounded}
(cf.\ \citet[\S7.2.B]{herlerLass99}).

It is known (cf.\ \citet[theorem 8]{heilmann77} and \citet[theorem 6.3.7]{herlerLass96}) that $\mathcal{M}_{\nu_0}$ in (\ref{eq:psop}) is a compact and convex region of $\mathbb{M}_w(\mathsf{K})$, which is spanned by stationary randomized policies and is characterized by linear constraints as
\begin{equation}
\label{eq:lcsaom}
\mathcal{M}_{\nu_0} =
\{\mu_{\nu_0}^{\pi}\colon \pi \in \Pi^{\scriptscriptstyle \textup{SR}}\} = 
\{\mu \in \mathbb{M}_w^{\scriptscriptstyle +}(\mathsf{K})\colon
\mathcal{L} \mu =  \nu_0\},
\end{equation}
where $\mathbb{M}_w^{\scriptscriptstyle +}(\mathsf{K}) \triangleq \{\mu  \in \mathbb{M}_w(\mathsf{K})\colon \mu \geqslant 0\}$ is the \emph{cone of measures} in $\mathbb{M}_w(\mathsf{K})$. 

We can thus explicitly formulate problem $L_{\lambda}(\nu_0)$ in (\ref{eq:psop}) as the infinite-dimensional
 LP problem
\begin{equation}
\label{eq:plp}
\begin{split}
L_{\lambda}(\nu_0)\colon \qquad & \maxim \; \langle \mu, r - \lambda c\rangle \\
& \st\colon \mathcal{L} \mu =  \nu_0, \mu \in \mathbb{M}_w^{\scriptscriptstyle +}(\mathsf{K}).
\end{split}
\end{equation}

\section{Decomposition of performance metrics.}
\label{s:sdr}
We present in this section a key decomposition of project performance metrics, which will be used in  subsequent analyses. 

We need a preliminary result. 
Note that the linear operator $\mathcal{L}^\star$ in (\ref{eq:tbetastu}) is (cf.\ \citet[p.\ 139]{herlerLass96}) the \emph{adjoint} of $\mathcal{L}$ in (\ref{eq:loL}), as it satisfies the \emph{basic adjoint relation}
\begin{equation}
\label{eq:ttstar}
\langle \mathcal{L} \mu, v\rangle = \langle \mu, \mathcal{L}^\star v\rangle, \quad 
\mu \in \mathbb{M}_w(\mathsf{K}), v \in  \mathbb{B}_w(\mathsf{X}),
\end{equation}
or,  in terms of the operators $\mathcal{L}^a$ in (\ref{eq:tbadef}),
\begin{equation}
\label{eq:tat}
\sum_{a \in \{0, 1\}} \langle \mathcal{L}^a \mu, v\rangle = 
\langle \mu, \mathcal{L}^\star v\rangle.
\end{equation}
The next result (which the author has not found in the literature) gives a
decomposition of  
(\ref{eq:tat}). 

\begin{lemma}
\label{lma:adjrel} 
For any $a \in
\{0, 1\}$, 
$\mu \in \mathbb{M}_w(\mathsf{K})$ and $v \in  \mathbb{B}_w(\mathsf{X}),$
\[
\langle \mathcal{L}^a \mu, v\rangle = 
\langle \tilde{\mu}^a, \mathcal{L}^\star v(\cdot, a)\rangle.
\] 
\end{lemma}
\proof
Let $a \in
\{0, 1\}$ and
$\mu \in \mathbb{M}_w(\mathsf{K})$.
For any $S \in \mathcal{B}(\mathsf{X})$ we have
\begin{align*}
\langle \mathcal{L}^a \mu, 1_S\rangle & = \mathcal{L}^a \mu(S) = \int \Big\{1_S(x) - \beta 
\int 1_S(y) \, \kappa^{a}(x, dy)\Big\} \, \tilde{\mu}^a(dx) \\
& = \int \mathcal{L}^\star 1_S (x, a) \,
\tilde{\mu}^a(dx) = \langle \tilde{\mu}^a, \mathcal{L}^\star 1_S (\cdot, a)\rangle,
\end{align*}
where we have used  (\ref{eq:bfnuv}), (\ref{eq:tbadef}),
(\ref{eq:tbetastu}), and (\ref{eq:bfmuu}). By standard arguments, it follows from the above that 
$\langle \mathcal{L}^a \mu, v\rangle = 
\langle \tilde{\mu}^a, \mathcal{L}^\star v (\cdot, a)\rangle$ for any 
$v \in \mathbb{B}_w(\mathsf{X})$.
 \endproof

We next give the main result of this section, consisting of two parts.
Part (a) decomposes
the reward metric $F(\nu_0, \pi)$ as the sum of $F(\nu_0, B)$
for a given active region $B$ 
  and 
a linear combination of  
$\tilde{\mu}_{\nu_0}^{\pi, a}$.
Part 
(b) similarly decomposes  the
resource metric $G(\nu_0, \pi)$.
This result extends to real-state projects corresponding  results for the discrete-state case in \citet[theorem 3]{nmaap01},
\citet[proposition 6.1]{nmmp02} and \citet[lemma 5.4]{nmmor06}.
Let $\nu_0 \in \mathbb{P}_w(\mathsf{X})$.

\begin{lemma}
\label{lma:dls}  For any policy $\pi \in \Pi$ and active region $B,$
\begin{itemize}
\item[\textup{(a)}]
$F(\nu_0, \pi)  = 
F(\nu_0, B) - \int_{B} f(y, B) \, \tilde{\mu}_{\nu_0}^{\pi, 0}(dy) + \int_{B^c} f(y, B) \, \tilde{\mu}_{\nu_0}^{\pi, 1}(dy),$
viz.,
\begin{equation}
\label{Fpidecref}
F(\nu_0, \pi) = F(\nu_0, B) +
\Ex_{\nu_0}^{\pi}\Bigg[\sum_{t=0}^\infty \beta^t \{A_t - 1_B(X_t)\} \, f(X_t, B)\Bigg];
\end{equation}
\item[\textup{(b)}]
$G(\nu_0, \pi)  = 
G(\nu_0, B) - \int_{B} g(y, B) \,
\tilde{\mu}_{\nu_0}^{\pi, 0}(dy) + \int_{B^c} g(y, B) \,
\tilde{\mu}_{\nu_0}^{\pi, 1}(dy),$
viz.,
\begin{equation}
\label{Gpidecref}
G(\nu_0, \pi) = G(\nu_0, B) +
\Ex_{\nu_0}^{\pi}\left[\sum_{t=0}^\infty \beta^t \{A_t - 1_B(X_t)\} \, g(X_t, B)\right].
\end{equation}
\end{itemize}
\end{lemma}
\proof
(a)
Using in turn Lemma \ref{lma:adjrel}, (\ref{eq:bfnuv}), 
(\ref{eq:fxBeq}) and
(\ref{eq:fseqs}), we can write 
\begin{align*}
\langle \mathcal{L}^0 \mu_{\nu_0}^\pi, F(\cdot, B)\rangle
 & = \left\langle \tilde{\mu}_{\nu_0}^{\pi, 0}, 
\mathcal{L}^\star F(\cdot, B)(\cdot, 0)\right\rangle
 = \int \mathcal{L}^\star F(\cdot, B)(y, 0) \, \tilde{\mu}_{\nu_0}^{\pi, 0}(dy) \\
& =  
\int_{B} \mathcal{L}^\star F(\cdot, B)(y, 0) \,
\tilde{\mu}_{\nu_0}^{\pi, 0}(dy) +
\int_{B^c} \mathcal{L}^\star F(\cdot, B)(y, 0) \,
 \tilde{\mu}_{\nu_0}^{\pi, 0}(dy)
\\
& = \int_{B} \{r(y, 0) + f(y, B)\} \,
\tilde{\mu}_{\nu_0}^{\pi, 0}(dy) + \int_{B^c} r(y, 0) \, \tilde{\mu}_{\nu_0}^{\pi, 0}(dy)
 \\
& = \int r(y, 0) \, \tilde{\mu}_{\nu_0}^{\pi, 0}(dy) 
+ \int_{B} f(y, B) \, \tilde{\mu}_{\nu_0}^{\pi, 0}(dy)
\end{align*}
and 
\begin{align*}
\langle \mathcal{L}^1 \mu_{\nu_0}^\pi, F(\cdot, B)\rangle &= \left\langle
\tilde{\mu}_{\nu_0}^{\pi, 1}, \mathcal{L}^\star F(\cdot, B)(\cdot, 1)\right\rangle = \int \mathcal{L}^\star F(\cdot, B)(y, 1) \,
 \tilde{\mu}_{\nu_0}^{\pi, 1}(dy) \\
& = \int_{B} \mathcal{L}^\star F(\cdot,
B)(y, 1) \, \tilde{\mu}_{\nu_0}^{\pi, 1}(dy) + 
\int_{B^c} \mathcal{L}^\star F(\cdot, B)(y, 1) \, 
\tilde{\mu}_{\nu_0}^{\pi, 1}(y)
 \\
& =  \int_{B}
r(y, 1) \, \tilde{\mu}_{\nu_0}^{\pi, 1}(dy) + \int_{B^c} \{r(y, 1) -
  f(y, B)\} \, \tilde{\mu}_{\nu_0}^{\pi, 1}(dy)  \\
& = \int  r(y, 1) \, \tilde{\mu}_{\nu_0}^{\pi, 1}(dy) -
\int_{B^c} f(y, B) \, \tilde{\mu}_{\nu_0}^{\pi, 1}(dy).
\end{align*}

Further, using that $\mathcal{L} = \sum_{a  \in \{0, 1\}} \mathcal{L}^a$, (\ref{eq:lcsaom}) and
(\ref{eq:bfnuv}), we
obtain
\begin{align*}
\sum_{a  \in \{0, 1\}} \langle \mathcal{L}^a \mu_{\nu_0}^\pi, F(\cdot, B)\rangle
& = \langle \mathcal{L} \mu_{\nu_0}^\pi, F(\cdot, B)\rangle = \langle \nu_0, F(\cdot, B)\rangle = F(\nu_0, B).
\end{align*}

Combining the above identities with 
\[
F(\nu_0, \pi) = \int  r \, d\mu_{\nu_0}^{\pi} = \sum_{a  \in \{0, 1\}} \int  r(y, a) \, \tilde{\mu}_{\nu_0}^{\pi, a}(dy).
\]
we readily obtain the stated decomposition identity
\begin{equation}
\label{eq:Fdec}
F(\nu_0, \pi)  = 
F(\nu_0, B) - \int_{B} f(y, B) \, \tilde{\mu}_{\nu_0}^{\pi, 0}(dy) + \int_{B^c} f(y, B) \, \tilde{\mu}_{\nu_0}^{\pi, 1}(dy).
\end{equation}

As for (\ref{Fpidecref}), it follows by reformulating (\ref{eq:Fdec}) 
using that 
\[
\int_{B} f(y, B) \, \tilde{\mu}_{\nu_0}^{\pi, 0}(dy) = \Ex_{\nu_0}^{\pi}\left[\sum_{t=0}^\infty \beta^t
  1_{B \times \{0\}}(X_t, A_t) \, f(X_t, B) \right]
\]
and
\[
 \int_{B^c} f(y, B) \, \tilde{\mu}_{\nu_0}^{\pi, 1}(dy) = \Ex_{\nu_0}^{\pi}\left[\sum_{t=0}^\infty \beta^t
  1_{B^c \times \{1\}} (X_t, A_t) \, f(X_t, B)\right]
\]
as
\[
F(\nu_0, \pi) +  \Ex_{\nu_0}^{\pi}\left[\sum_{t=0}^\infty \beta^t
  1_{B \times \{0\}}(X_t, A_t) \, f(X_t, B) \right] 
  = F(\nu_0, B) + \Ex_{\nu_0}^{\pi}\left[\sum_{t=0}^\infty \beta^t
  1_{B^c \times \{1\}} (X_t, A_t) \, f(X_t, B)\right],
\]
and then reformulating the latter identity as (\ref{Fpidecref}) noting that, for any state $x$ and action $a$,
\[
  1_{B^c \times \{1\}} (x, a)  - 
  1_{B \times \{0\}}(x, a)  = 
  a - 1_B(x).
\]

Part (b) follows along the same lines as part (a).
 \endproof

We will use Lemma \ref{lma:dls} as a powerful tool to generate useful identities.

\section{Analysis of metrics as functions of the threshold.}
\label{s:ppmftv} 
This section presents properties of metrics under threshold policies as functions of the threshold.
Throughout, $\nu_0 \in \mathbb{P}_w(\mathsf{X})$ and $x \in \mathsf{X}$.
\subsection{C\`adl\`ag property.}
\label{s:cadlagp}
Lemma \ref{pro:cadlagFG} shows that  $F(\nu_0, z),$ $G(\nu_0, z),$ $f(x, z)$ and $g(x, z)$ are c\`adl\`ag in $z \in \mathbb{R}$.
We need the following result, referring to 
the extension to $\mathbb{R}$ of the state-occupation measures
$\tilde{\mu}_{\nu_0}^{\pi}$, viz., $\tilde{\mu}_{\nu_0}^{\pi}(S) \triangleq \tilde{\mu}_{\nu_0}^{\pi}(S \cap \mathsf{X})$ for $S \in \mathcal{B}(\mathbb{R})$. 
We abuse notation by letting, e.g., $[z_1,  z_2]$ mean 
$\{z \in \mathbb{R}\colon z_1 \leqslant z \leqslant z_2\}$, even when $z_1 = -\infty$ or $z_2 = \infty$, and similarly with other intervals.

\begin{lemma}
\label{lma:FGxzdfgxz} For any thresholds $-\infty \leqslant z_1 < z_2 \leqslant \infty,$ 
\begin{itemize}
\item[\textup{(a)}] $F(\nu_0, z_1) -  F(\nu_0, z_2) = \int_{(z_1,  z_2]} f(y, z_2) \, \tilde{\mu}_{\nu_0}^{z_1}(dy) = \int_{(z_1, z_2]} f(y, z_1) \, \tilde{\mu}_{\nu_0}^{z_2}(dy);$
\item[\textup{(b)}] $G(\nu_0, z_1) -  G(\nu_0, z_2) = \int_{(z_1, z_2]}
  g(y, z_2) \, \tilde{\mu}_{\nu_0}^{z_1}(dy) = \int_{(z_1, z_2]}
  g(y, z_1) \, \tilde{\mu}_{\nu_0}^{z_2}(dy);$
\item[\textup{(c)}] $F(\nu_0, z_1) -
  F(\nu_0, z_2^{\scriptscriptstyle -}) = \int_{(z_1, z_2)} f(y, z_2^{\scriptscriptstyle -}) \,
  \tilde{\mu}_{\nu_0}^{z_1}(dy) = \int_{(z_1, z_2)} f(y, z_1) \,
  \tilde{\mu}_{\nu_0}^{z_2^{\scriptscriptstyle -}}(dy);$
\item[\textup{(d)}] $G(\nu_0, z_1) -
  G(\nu_0, z_2^{\scriptscriptstyle -}) =  \int_{(z_1, z_2)} g(y, z_2^{\scriptscriptstyle -}) \,
  \tilde{\mu}_{\nu_0}^{z_1}(dy) = \int_{(z_1, z_2)} g(y, z_1) \,
  \tilde{\mu}_{\nu_0}^{z_2^{\scriptscriptstyle -}}(dy).$
  \item[\textup{(e)}] $F(\nu_0, z_1^{\scriptscriptstyle -}) -  F(\nu_0, z_2) = \int_{[z_1, z_2]} f(y, z_2) \, \tilde{\mu}_{\nu_0}^{z_1^{\scriptscriptstyle -}}(dy) = \int_{[z_1, z_2]} f(y, z_1^{\scriptscriptstyle -}) \, \tilde{\mu}_{\nu_0}^{z_2}(dy);$
  \item[\textup{(f)}] $G(\nu_0, z_1^{\scriptscriptstyle -}) -
  G(\nu_0, z_2) =  \int_{[z_1, z_2]} g(y, z_2^{\scriptscriptstyle -}) \,
  \tilde{\mu}_{\nu_0}^{z_1^{\scriptscriptstyle -}}(dy) = \int_{[z_1, z_2]} g(y, z_1^{\scriptscriptstyle -}) \,
  \tilde{\mu}_{\nu_0}^{z_2}(dy).$
\end{itemize}
\end{lemma}
\proof
(a, b) For the first identity take $\pi = z_1$ and $B = (z_2, u]$ in Lemma
\ref{lma:dls}(a, b), using that $\tilde{\mu}_{\nu_0}^{z_1, 0}(z_1, u] = \tilde{\mu}_{\nu_0}^{z_1, 1}[\ell, z_1] = 0$ and $\tilde{\mu}_{\nu_0}^{z_1,
  1} = \tilde{\mu}_{\nu_0}^{z_1}$ on $\mathcal{B}(z_1, z_2]$. For the second identity take $\pi = z_2$ and $B = (z_1, u]$, using that $\tilde{\mu}_{\nu_0}^{z_2, 0}(z_2, u] = \tilde{\mu}_{\nu_0}^{z_2, 1}[\ell, z_1] = 0$ and $\tilde{\mu}_{\nu_0}^{z_2,
  0} = \tilde{\mu}_{\nu_0}^{z_2}$ on $\mathcal{B}(z_1, z_2]$. 

(c, d) For the first identity take $\pi =  z_1$ and $B = [z_2, u]$ in Lemma
\ref{lma:dls}(a, b), using that $\tilde{\mu}_{\nu_0}^{z_1, 0}[z_2, u] = \tilde{\mu}_{\nu_0}^{z_1, 1}[\ell, z_1] = 0$ and $\tilde{\mu}_{\nu_0}^{z_1,
  1} = \tilde{\mu}_{\nu_0}^{z_1}$ on $\mathcal{B}(z_1, z_2)$.
 For the second identity take $\pi = z_2^{-}$ and $B = (z_1, u]$, using that
$\tilde{\mu}_{\nu_0}^{z_2^{\scriptscriptstyle -}, 0}[z_2, u] =
\tilde{\mu}_{\nu_0}^{z_2^{\scriptscriptstyle -}, 1}[\ell, z_1] = 0$
and $\tilde{\mu}_{\nu_0}^{z_2^{\scriptscriptstyle -}, 0} =
\tilde{\mu}_{\nu_0}^{z_2^{\scriptscriptstyle -}}$ on
$\mathcal{B}(z_1, z_2)$. 

(e, f) For the first identity take $\pi =  z_1^{-}$ and $B = (z_2, u]$ in Lemma
\ref{lma:dls}(a, b), using that $\tilde{\mu}_{\nu_0}^{z_1^{\scriptscriptstyle -}, 0}(z_2, u] = \tilde{\mu}_{\nu_0}^{z_1^{\scriptscriptstyle -}, 1}(\bar{\mathsf{X}}_{z_1}^c) = 0$ and $\tilde{\mu}_{\nu_0}^{z_1^{\scriptscriptstyle -},
  1} = \tilde{\mu}_{\nu_0}^{z_1^{\scriptscriptstyle -}}$ on $\mathcal{B}[z_1, z_2]$.
 For the second identity take $\pi = z_2$ and $B = [z_1, u]$, using that
$\tilde{\mu}_{\nu_0}^{z_2, 0}(z_2, u] =
\tilde{\mu}_{\nu_0}^{z_2, 1}[\ell, z_1) = 0$
and $\tilde{\mu}_{\nu_0}^{z_2, 0} =
\tilde{\mu}_{\nu_0}^{z_2}$ on
$\mathcal{B}[z_1, z_2]$. 
 \qquad
 \endproof

Recall from Section \ref{s:assub} that the notation $F(\nu_0, z^{\scriptscriptstyle -})$ (for this and other metrics) refers to $F(\nu_0, \pi)$ under the $z^{\scriptscriptstyle -}$-policy.
The next result ensures the consistency of  such notation with the standard one where $F(\nu_0, z^{\scriptscriptstyle -})$ denotes the left limit at $z$ of function $F(\nu_0, \cdot)$.
The result further shows the consistency of the notation $F(\nu_0, -\infty)$ and $F(\nu_0, \infty)$, where the second argument refers to the threshold policies $-\infty$ (``always active'') and $\infty$ (``never active'') with the standard notation where $F(\nu_0, -\infty)$ and $F(\nu_0, \infty)$ are the limits of $F(\nu_0, z)$ as $z \to -\infty$ and as $z \to \infty$, respectively.

\begin{lemma}
\label{pro:cadlagFG} The following holds$:$
\begin{itemize}
\item[\textup{(a)}]
$F(\nu_0, \cdot),$ $G(\nu_0, \cdot),$ $f(x, \cdot)$ and $g(x, \cdot)$ are bounded c\`adl\`ag functions with
left limits at $z \in \mathbb{R}$ given by $F(\nu_0, z^{\scriptscriptstyle -}),$ $G(\nu_0, z^{\scriptscriptstyle -}),$ $f(x, z^{\scriptscriptstyle -})$ and $g(x,
z^{\scriptscriptstyle -});$
\item[\textup{(b)}] $\lim_{z \to -\infty} \, (F(\nu_0, z), G(\nu_0, z), f(x, z), g(x, z))  = (F(\nu_0, -\infty), G(\nu_0, -\infty), f(x, -\infty), g(x, -\infty));$
\item[\textup{(c)}] $\lim_{z \to \infty} \, (F(\nu_0, z), G(\nu_0, z), f(x, z), g(x, z))  = (F(\nu_0, \infty), G(\nu_0, \infty), f(x, \infty), g(x, \infty)).$
\end{itemize}
\end{lemma}
\proof (a) We first show that $F(\nu_0, \cdot)$ is right-continuous at $z$.
From Lemma \ref{lma:FGxzdfgxz}(a) we obtain (with $M_\gamma$ as in (\ref{eq:FGMwnon})), for $\delta > 0$, 
\begin{align*}
\big|F(\nu_0, z+\delta) - F(\nu_0, z)\big| & \leqslant  \int_{(z, z+\delta]} |f|(y, z+\delta) \, \tilde{\mu}_{\nu_0}^{z}(dy)   \leqslant 2 M_\gamma \int_{(z, z+\delta]} w(y) \, \tilde{\mu}_{\nu_0}^{z}(dy) \to 0 \textup{ as } \delta \searrow 0,
\end{align*}
where
the second inequality follows from (\ref{eq:fgxpibnd}), and the
limit  follows by the dominated convergence
theorem, since $\tilde{\mu}_{\nu_0}^{z}  \in \mathbb{M}_w(\mathsf{X})$ (see Section \ref{s:saom}).
The result for $G(\nu_0, \cdot)$ follows similarly. Furthermore, 
\begin{align*}
\lim_{\delta \searrow 0} f(x, z+\delta) & = 
\lim_{\delta \searrow 0}  \diff\limits_{a = 0}^{a=1} \Big\{r(x, a) + \beta \int F(y, z+\delta) \, \kappa^a(x, dy) \Big\} \\
& = \diff\limits_{a = 0}^{a=1}  \Big\{r(x, a) + \beta \int F(y, z) \, \kappa^a(x, dy) \Big\}  = f(x, z),
\end{align*}
where the interchange of limit and integral is justified by the
dominated convergence theorem, using  (\ref{eq:FGMwnon}) and the right-continuity of $F(y, \cdot)$.
The result for $g(x, \cdot)$ follows similarly. 

Consider now the left limits. 
From Lemma \ref{lma:FGxzdfgxz}(c) and arguing as above we obtain
\begin{align*}
\big|F(\nu_0, z-\delta) - F(\nu_0, z^{\scriptscriptstyle -})\big| & \leqslant  \int_{(z-\delta, z)} |f|(y, z-\delta) \, \tilde{\mu}_{\nu_0}^{z^{\scriptscriptstyle -}}(dy)   \leqslant 2 M_\gamma \int_{(z-\delta, z)} w(y) \, \tilde{\mu}_{\nu_0}^{z^{\scriptscriptstyle -}}(dy) \to 0 \textup{ as } \delta \searrow 0,
\end{align*}
and similarly for the left limits of $G(\nu_0, \cdot)$,  
$f(x, \cdot)$ and 
$g(x, \cdot)$ at $z$.

(b) The result is trivial if $\ell$ is finite, so suppose $\ell = -\infty.$ For any finite $z$, taking $z_1 = -\infty$ and $z_2 = z$ in the first identity in Lemma \ref{lma:FGxzdfgxz}(a) yields
\begin{equation}
\label{eq:Fxminfz}
F(\nu_0, -\infty) = F(\nu_0, z) + \int_{(-\infty, z]} f(y, z) \, \tilde{\mu}_{\nu_0}^{-\infty}(dy).
\end{equation}
The result that $\lim_{z \to -\infty} F(\nu_0, z) = F(\nu_0, -\infty)$ now follows from (\ref{eq:Fxminfz}) by the dominated convergence theorem, arguing as in part (a).
The results for $G(\nu_0, -\infty)$, $f(x, -\infty)$ and $g(x, -\infty)$ follow similarly.

(c) The result is trivial if $u$ is finite, so suppose $u = \infty.$ For any finite $z$, taking $z_1 = z$ and $z_2 = \infty$ in the second identity in Lemma \ref{lma:FGxzdfgxz}(a) yields
\begin{equation}
\label{eq:Fxminfz2}
F(\nu_0, \infty) = F(\nu_0, z) - \int_{(z, \infty)} f(y, z) \, \tilde{\mu}_{\nu_0}^{\infty}(dy).
\end{equation}
The result that $\lim_{z \to \infty} F(\nu_0, z) = F(\nu_0, \infty)$ follows from (\ref{eq:Fxminfz2}) by the dominated convergence theorem, arguing as in part (a).
The results for $G(\nu_0, \infty)$, $f(x, \infty)$ and $g(x, \infty)$ follow similarly.
\qquad
 \endproof

\subsection{Monotonicity of resource metric.}
\label{s:mbbp}
The next result shows that, under (PCLI1), the resource metric $G(\nu_0, z)$ is, as intuition would suggest, monotone nonincreasing in the threshold $z$. 

\begin{lemma}
\label{pro:Gnondecr}
Let \textup{(PCLI1)} hold$.$ Then
\begin{itemize}
\item[\textup{(a)}] $G(\nu_0, z)$ is nonincreasing  in $z$ on $\mathbb{R};$
\item[\textup{(b)}] if $\nu_0$ has full support $\mathsf{X}$ then $G(\nu_0, z)$ is   decreasing in $z$ on $\mathsf{X}.$
\end{itemize}
\end{lemma}
\proof
(a) For finite $z_1 < z_2$, 
 Lemma
\ref{lma:FGxzdfgxz}(b) and  (PCLI1) yield $G(\nu_0, z_1) \geqslant  G(\nu_0,
z_2)$.

(b) 
Lemma
\ref{lma:FGxzdfgxz}(b) further implies that, for any finite thresholds $z_1 < z_2$ in $\mathsf{X}$, 
\begin{equation}
\label{eq:Gnu0zpzeq}
G(\nu_0, z_1) -  G(\nu_0, z_2) = \int_{(z_1, z_2]}
  g(y, z_2) \, \tilde{\mu}_{\nu_0}^{z_1}(dy) > 0,
\end{equation}
where the inequality follows from 
 (PCLI1) and $\tilde{\mu}_{\nu_0}^{z_1}(z_1, z_2] \geqslant 
\nu_0(z_1, z_2]> 0$.
\qquad 
 \endproof

We next establish 
 existence and finiteness of the LS integrals in condition (PCLI3).

\begin{lemma}
\label{lma:eaafpcli3} Under \textup{(PCLI1, PCLI2)}$,$ $\int_{(z_1, z_2]} m(z) \, G(\nu_0, dz)$ exists and is finite for finite $z_1 < z_2.$ 
\end{lemma}
\proof
By
Lemmas \ref{pro:cadlagFG} and \ref{pro:Gnondecr}(a), $G(\nu_0, \cdot)$ is bounded, c\`adl\`ag  and nonincreasing. Since  $m(\cdot)$ is continuous, $\int_{(z_1, z_2]} m(z) \, G(\nu_0, dz)$ exists and is finite by standard  results on LS integration. 
 \endproof

\subsection{Analysis of discontinuities.}
\label{s:aod}
Being c\`adl\`ag, functions 
 $F(\nu_0, \cdot)$, $G(\nu_0, \cdot)$, $f(x, \cdot)$ and $g(x,
\cdot)$ have at most
countably many discontinuities, which are of jump type.
We analyze below such discontinuities.  The following result  
 gives simple formulae for their  jumps.
 We denote by $\Delta_{2} F(\nu_0, z) \triangleq F(\nu_0, z) - F(\nu_0, z^{\scriptscriptstyle -})$ the jump of  $F(\nu_0, \cdot)$ at $z$, and similarly for $\Delta_{2} G(\nu_0, z)$, $\Delta_{2} f(x, z)$ and $\Delta_{2} g(x, z)$.

\begin{lemma} 
\label{lms:gfminus} For any threshold $z \in \mathsf{X},$ 
\begin{itemize}
\item[\textup{(a)}] $\Delta_{2} F(\nu_0, z)  = -f(z,
  z^{\scriptscriptstyle -}) \, \tilde{\mu}_{\nu_0}^{z}\{z\} = -f(z, z)
  \, \tilde{\mu}_{\nu_0}^{z^{\scriptscriptstyle -}}\{z\};$
\item[\textup{(b)}] 
$\Delta_{2} G(\nu_0, z)  = -g(z, z^{\scriptscriptstyle -}) \, \tilde{\mu}_{\nu_0}^{z}\{z\} = -g(z, z)
  \, \tilde{\mu}_{\nu_0}^{z^{\scriptscriptstyle -}}\{z\}$.
\end{itemize}
\end{lemma}
\proof
(a) Taking $\pi = z$ and $B = \{y \in \mathsf{X}\colon y \geqslant z\}$ in Lemma
\ref{lma:dls}(a) gives
\[
F(\nu_0, z) + \int_{y \geqslant z} f(y, z^{\scriptscriptstyle -}) \, \tilde{\mu}_{\nu_0}^{z, 0}(dy) = F(\nu_0, z^{\scriptscriptstyle -}) + \int_{y < z} f(y, z^{\scriptscriptstyle -}) \, \tilde{\mu}_{\nu_0}^{z, 1}(dy),
\]
so $\Delta_{2} F(\nu_0, z)  + f(z, z^{\scriptscriptstyle -}) \tilde{\mu}_{\nu_0}^{z}\{z\} = 0$, as $\tilde{\mu}_{\nu_0}^{z, 0}\{z\} = \tilde{\mu}_{\nu_0}^{z}\{z\}$ and $\tilde{\mu}_{\nu_0}^{z, 0}\{y >z\} = \tilde{\mu}_{\nu_0}^{z, 1}\{y < z\} = 0$. 

Further, taking $\pi = z^{\scriptscriptstyle -}$ and $B = \{y \in \mathsf{X}\colon y > z\}$ in Lemma
\ref{lma:dls}(a) gives
\[
F(\nu_0, z^{\scriptscriptstyle -}) + \int_{y > z} f(y, z) \, \tilde{\mu}_{\nu_0}^{z^{\scriptscriptstyle -}, 0}(dy) = F(\nu_0, z) + \int_{y \leqslant z} f(y, z) \, \tilde{\mu}_{\nu_0}^{z^{\scriptscriptstyle -}, 1}(dy),
\]
so  $\Delta_{2} F(\nu_0, z)  + f(z, z) \, \tilde{\mu}_{\nu_0}^{z^{\scriptscriptstyle -}}\{z\} = 0$, as $\tilde{\mu}_{\nu_0}^{z^{\scriptscriptstyle -}, 1}\{z\} = \tilde{\mu}_{\nu_0}^{z^{\scriptscriptstyle -}}\{z\}$ and $\tilde{\mu}_{\nu_0}^{z^{\scriptscriptstyle -}, 0}\{y > z\} = \tilde{\mu}_{\nu_0}^{z^{\scriptscriptstyle -}, 1}\{y < z\} = 0$. 

Part (b) follows along the same lines as part (a).
 \endproof

Recall that PCL-indexability condition \textup{(PCLI1)} requires $g(x, z)$ to be positive for any state $x$ and threshold
$z$.
The following result shows that it further implies $g(x, x^{\scriptscriptstyle -}) > 0$.

\begin{lemma}
\label{lma:gzzmp}
Under \textup{(PCLI1),} 
$g(x, x^{\scriptscriptstyle -}) > 0$ and $\Delta_{2} G(x, x) < 0.$
\end{lemma}
\proof
Taking $\nu_0 = \delta_x$ and $z = x$ in Lemma \ref{lms:gfminus}(b) and using  (PCLI1) and 
$\tilde{\mu}_{x}^{x^{\scriptscriptstyle -}}\{x\} \geqslant 1$
gives 
\[
g(x, x^{\scriptscriptstyle -}) \, \tilde{\mu}_{x}^{x}\{x\} = g(x, x)
  \, \tilde{\mu}_{x}^{x^{\scriptscriptstyle -}}\{x\} \geqslant g(x, x) > 0,
\] 
whence $g(x, x^{\scriptscriptstyle -}) > 0$ and $\Delta_{2} G(x, x) < 0$. \qquad 
 \endproof

The next result characterizes the discontinuities of
$G(\nu_0, \cdot)$ in terms of $T_z \triangleq \min \{t \geqslant 0\colon X_t = z\}$, the \emph{first hitting time to $z$} of state process $\{X_t\}_{t=0}^\infty$ under the $z$-policy,
with 
$T_z \triangleq \infty$ if $z$ is never hit.

\begin{lemma}
\label{pro:gfcont}
Under \textup{(PCLI1)}$,$ 
$G(\nu_0, \cdot)$ is discontinuous at $z$ iff $\Prob_{\nu_0}^{z}\{T_{z} < \infty\} > 0.$
\end{lemma}
\proof
It follows from Lemmas \ref{lms:gfminus}(b) and 
\ref{lma:gzzmp} that $G(\nu_0, \cdot)$ is discontinuous at $z$
iff $\tilde{\mu}_{\nu_0}^{z}\{z\} > 0$, which happens iff
$\Prob_{\nu_0}^{z}\{T_{z} < \infty\} > 0$.
 \endproof

The next result relates the  jumps of $F(\nu_0, \cdot)$ and $G(\nu_0, \cdot)$, gives a corresponding result for $f(x, \cdot)$ and $g(x, \cdot)$, and ensures that the MP metric $m(x, \cdot)$ (see (\ref{eq:mpm1}))  is continuous at $z = x$.

\begin{lemma}
\label{pro:mafgids} 
Let \textup{(PCLI1)} hold$.$ Then$,$ for any threshold $z \in \mathsf{X},$ 
\begin{itemize}
\item[\textup{(a)}]  $\Delta_{2} F(\nu_0, z)  =
  m(z) \, \Delta_{2} G(\nu_0, z);$
\item[\textup{(b)}] $\Delta_{2} f(x, z)  = m(z)
  \, \Delta_{2} g(x, z);$
\item[\textup{(c)}] $m(x, \cdot)$ is continuous at $x,$ viz$.,$ $m(x, x^{\scriptscriptstyle -}) = m(x, x) = m(x).$
\end{itemize}
\end{lemma}
\proof
(a)  Lemma \ref{lms:gfminus}(a, b) and $g(z, z) > 0$ yield
\begin{align*}
\Delta_{2} F(\nu_0, z)  = -f(z, z)
\, \tilde{\mu}_{\nu_0}^{z^{\scriptscriptstyle -}}\{z\} = 
-m(z) \, g(z, z) \, \mu_{\nu_0}^{z^{\scriptscriptstyle -}}\{z\} = 
m(z) \Delta_{2} G(\nu_0, z) .
\end{align*}

(b) Using in turn (\ref{eq:mrm1}), part (a) and (\ref{eq:mwm1}), gives
\begin{align*}
\Delta_{2} f(x, z)   & = \beta \diff\limits_{a=0}^{a=1} \int
\Delta_{2} F(y, z) \,
\kappa^{a}(x, dy)
 = 
m(z) \beta \diff\limits_{a=0}^{a=1} \int  \Delta_{2} G(y, z)
\kappa^{a}(x, dy) = m(z) \Delta_{2} g(x, z). 
\end{align*}

(c) Taking $z = x$ in part (b) and noting that $m(x) = f(x,
x)/g(x, x)$, we obtain
\[
f(x, x^{\scriptscriptstyle -}) - f(x, x) = m(x) [g(x, x^{\scriptscriptstyle -}) -
  g(x, x)] = m(x) g(x, x^{\scriptscriptstyle -}) - f(x, x).
\]
whence $f(x, x^{\scriptscriptstyle -})=
m(x) g(x, x^{\scriptscriptstyle
  -})$ and, by Lemma \ref{lma:gzzmp}, 
$m(x) =
m(x, x^{\scriptscriptstyle
  -})$.
\qquad
 \endproof

Lemma \ref{pro:mafgids}(a) yields
 an alternate characterization of the MP index under (PCLI1) in terms of jumps $\Delta_{2} F(x, x) = F(x, x) - F(x, x^{\scriptscriptstyle -})$ and $\Delta_{2} G(x, x) = G(x, x) - G(x, x^{\scriptscriptstyle -})$ (cf.\ Lemma \ref{lma:gzzmp}).

\begin{corollary}
\label{cor:mpiac}  
 Under \textup{(PCLI1)}$,$ 
$m(x) = \Delta_{2} F(x, x)/\Delta_{2} G(x, x).$
\end{corollary}

Further, Lemma \ref{pro:mafgids}(a, b) immediately implies  the following result relating the  continuity points of functions $F(x, \cdot)$ and $G(x, \cdot)$, and of  $f(x, \cdot)$ and $g(x, \cdot)$.

\begin{corollary}
\label{cor:fgids} 
Under \textup{(PCLI1),} 
\begin{itemize}
\item[\textup{(a)}] if $G(\nu_0, \cdot)$ is continuous at $z,$ so is $F(\nu_0, \cdot);$ 
\item[\textup{(b)}] if $g(x, \cdot)$ is continuous at $z,$ so is $f(x, \cdot).$ 
\end{itemize}
\end{corollary}

\subsection{Bounded variation and MP index as Radon--Nikodym derivative.}
\label{s:bvp}
We next show that the metrics of a PCL-indexable project, as functions of the  threshold, belong to the linear space  $\mathbb{V}(\mathbb{R})$ of functions of \emph{bounded variation}  on $\mathbb{R}$ (see \citet[\S2.7]{cartvanBrunt00}). 
We further build on that result to characterize the MP index as a Radon--Nikodym derivative.

We start with the resource metrics, for which condition (PCLI1) suffices. 

\begin{lemma}
\label{pro:gbv} 
Under \textup{(PCLI1),} $G(\nu_0, \cdot), g(x, \cdot) \in \mathbb{V}(\mathbb{R}).$ 
\end{lemma}
\proof 
$G(\nu_0, \cdot) \in \mathbb{V}(\mathbb{R})$ follows from
$G(\nu_0, \cdot)$ being bounded  (see (\ref{eq:FGMwnon})) and nonincreasing on $\mathbb{R}$ (see Lemma \ref{pro:Gnondecr}(a)).
Since
$
g(x, z) = \diff_{a=0}^{a=1} [c(x, a) + \beta G(\kappa^a(x, \cdot), z)],
$
$g(x, \cdot)$ is a difference of bounded nonincreasing functions, whence (see \citet[theorem 2.7.2]{cartvanBrunt00}) $g(x, \cdot) \in \mathbb{V}(\mathbb{R}).$
 \endproof

\begin{remark}
\label{re:nuGx}
Since the function $G(\nu_0, \cdot)$ is
bounded, c\`adl\`ag and nonincreasing$,$ by Carath\'eodory's extension theorem
it induces a unique \emph{finite LS measure} $\nu_{\scriptscriptstyle G(\nu_0, \cdot)}$ on $\mathcal{B}(\mathbb{R})$ satisfying
\[
\nu_{\scriptscriptstyle G(\nu_0, \cdot)}(z_1, z_2] = G(\nu_0, z_1) - G(\nu_0, z_2), \quad 
-\infty < z_1 < z_2 < \infty.
\]

\end{remark}

In light of Remark \ref{re:nuGx}, we have by standard results that the LS integral 
\[
\int |m|(z) \,G(\nu_0, dz) = -\int |m| \,d\nu_{\scriptscriptstyle G(\nu_0, \cdot)}
\] is well defined.
The following result shows that such an integral is finite, even when $\mathsf{X}$ is
 unbounded.

\begin{lemma}
\label{lma:phimpgint}
Under PCL-indexability$,$ $m(\cdot)$ is $\nu_{\scriptscriptstyle G(\nu_0, \cdot)}$-integrable$.$ 
\end{lemma}
\proof
Consider first the case that $m(\cdot) \geqslant 0$ over $\mathbb{R}$. Then,
\begin{align*}
\int |m| \, d\nu_{\scriptscriptstyle G(\nu_0, \cdot)} & =  \int m \, \nu_{\scriptscriptstyle G(\nu_0, \cdot)} = 
\lim_{n \to \infty} \, \int_{(-n, n]} m \, \nu_{\scriptscriptstyle G(\nu_0, \cdot)} \\ & = 
\lim_{n \to \infty} \, \{F(\nu_0, -n) - F(\nu_0, n)\} =
F(\nu_0, -\infty) - F(\nu_0, \infty)< \infty,
\end{align*}
where we have used the monotone convergence theorem, condition (PCLI3) and Lemma \ref{pro:cadlagFG}(b, c).
A similar argument yields the result in the case that $m \leqslant 0$ over $\mathbb{R}$.

Otherwise, 
 (PCLI2) implies that
there is a state $b$ with $m(b) = 0$ such that $|m| = - m$ on $(-\infty, b]$ and $|m| = m$ on $(b, \infty)$.
Arguing along the same lines as above yields
\begin{align*}
\int |m| \, \nu_{\scriptscriptstyle G(\nu_0, \cdot)} & = \int_{(-\infty, b]} -m \, d\nu_{\scriptscriptstyle G(\nu_0, \cdot)} + \int_{(b, \infty)} m \, \nu_{\scriptscriptstyle G(\nu_0, \cdot)}  = \lim_{n \to \infty} \, \int_{(-n, b]} -m \, d\nu_{\scriptscriptstyle G(\nu_0, \cdot)} + \lim_{n \to \infty} \, \int_{(b, n]} m \, \nu_{\scriptscriptstyle G(\nu_0, \cdot)} \\
& = \lim_{n \to \infty} \, \{F(\nu_0, b) - F(\nu_0, -n)\} + \lim_{n \to \infty} \, \{F(\nu_0, b) - F(\nu_0, n)\} \\
& = 2 F(\nu_0, b) - F(\nu_0, -\infty) - F(\nu_0, \infty) < \infty.
\end{align*}

\endproof

The following result shows that, under PCL-indexability, the identity (\ref{eq:pcli3nu}) ---and, in particular condition (PCLI3)--- extends to infinite intervals. 

\begin{lemma}
\label{lma:pcli3inf} Under PCL-indexability$,$ 
\begin{equation}
\label{eq:pcli3inf}
F(\nu_0, z_2) -  F(\nu_0, z_1) = \int_{z_1 < z \leqslant z_2} m(z) \, G(\nu_0, dz), \quad -\infty \leqslant z_1 < z_2   \leqslant \infty.
\end{equation}
\end{lemma}
\proof
The result follows from (\ref{eq:pcli3nu}), Lemmas \ref{pro:cadlagFG}(b) and \ref{lma:phimpgint}, and  dominated convergence.
 \endproof

We next present the counterpart of Lemma \ref{pro:gbv} for reward measures.

\begin{lemma}
\label{pro:fbv} Under PCL-indexability$,$ $F(\nu_0, \cdot), f(x, \cdot) \in \mathbb{V}(\mathbb{R}).$ 
\end{lemma}
\proof
Since $\nu_{\scriptscriptstyle G(\nu_0, \cdot)}$ (see Remark \ref{re:nuGx}) is a finite measure,  
Lemma \ref{lma:phimpgint} ensures
  that $\tilde{\nu}_{\nu_0}(S) \triangleq \int_{S} m \, d\nu_{\scriptscriptstyle G(\nu_0, \cdot)}$  is a finite signed measure,
 and hence admits a Jordan decomposition $\tilde{\nu}_{\nu_0} = \tilde{\nu}_{\nu_0}^{\scriptscriptstyle +} - \tilde{\nu}_{\nu_0}^{\scriptscriptstyle -}$ with $\tilde{\nu}_{\nu_0}^{\scriptscriptstyle +}$ and $\tilde{\nu}_{\nu_0}^{\scriptscriptstyle -}$  finite measures. 
On the other hand,  (PCLI3) and Lemma \ref{pro:cadlagFG}(b) yield $\tilde{\nu}_{\nu_0}(-\infty, z] =  F(\nu_0, -\infty) - F(\nu_0, z)$ for $z \in \mathbb{R}$,
whence
 $F(\nu_0, z)  = F(\nu_0, -\infty) + \tilde{\nu}_{\nu_0}^{\scriptscriptstyle -}(-\infty, z] - \tilde{\nu}_{\nu_0}^{\scriptscriptstyle +}(-\infty, z]$. Hence, $F(\nu_0, \cdot) \in \mathbb{V}(\mathbb{R})$, being a difference of  bounded nondecreasing functions. 
 Further,
it follows that function $z \mapsto F(\kappa^a(\nu_0, \cdot), z)$ is in $\mathbb{V}(\mathbb{R})$.
Thus, the identity
$
f(x, z) = \diff_{a=0}^{a=1} [r(x, a) + \beta F(\kappa^a(x, \cdot), z)]
$
 represents $f(x, \cdot)$ as the difference of two  functions in $\mathbb{V}(\mathbb{R})$, whence $f(x, \cdot) \in \mathbb{V}(\mathbb{R}).$
 \endproof

\begin{remark}
\label{re:nuFx}
In light of Lemma \ref{pro:fbv} and $F(\nu_0, \cdot)$ being bounded c\`adl\`ag  under PCL-indexability (see (\ref{eq:FGMnuwnon}) and Lemma \ref{pro:cadlagFG}), 
Carath\'eodory's theorem for signed measures (cf.\ \citet[section X.6]{doob94}) ensures  existence of a unique 
 \emph{finite LS signed measure} $\nu_{\scriptscriptstyle F(\nu_0, \cdot)}$ on $\mathcal{B}(\mathbb{R})$ satisfying
\[
\nu_{\scriptscriptstyle F(\nu_0, \cdot)}(z_1, z_2] \triangleq F(\nu_0, z_1) - F(\nu_0, z_2), \quad 
-\infty < z_1 < z_2 < \infty.
\]
 \end{remark}

The next result shows that the MP index is a Radon--Nikodym derivative of $\nu_{\scriptscriptstyle F(\nu_0, \cdot)}$ with respect to $\nu_{\scriptscriptstyle G(\nu_0, \cdot)}$. The notation $\ll$ below means ``is absolutely continuous with respect to''.

\begin{proposition}
\label{pro:mpirnd} 
Under PCL-indexability$,$ $\nu_{\scriptscriptstyle F(\nu_0, \cdot)} \ll \nu_{\scriptscriptstyle G(\nu_0, \cdot)},$
with $m(\cdot)$ being a Radon--Nikodym derivative of $\nu_{\scriptscriptstyle F(\nu_0, \cdot)}$ with respect to $\nu_{\scriptscriptstyle G(\nu_0, \cdot)},$ viz.$,$
\begin{equation}
\label{eq:ciii}
 \nu_{\scriptscriptstyle F(\nu_0, \cdot)}(S) = \int_{S} m(x) \, \nu_{\scriptscriptstyle G(\nu_0, \cdot)}(dx), \quad S \in \mathcal{B}(\mathbb{R}). 
\end{equation}
\end{proposition}
\proof
The result that $\nu_{\scriptscriptstyle F(\nu_0, \cdot)}$  satisfies (\ref{eq:ciii}) follows from (\ref{eq:pcli3}) and Remark \ref{re:nuFx}, whereas $\nu_{\scriptscriptstyle F(\nu_0, \cdot)} \ll \nu_{\scriptscriptstyle G(\nu_0, \cdot)}$ follows directly from (\ref{eq:ciii}). 
 \endproof

By standard results on differentiation of LS measures (cf.\ \citet[section X.4]{doob94}), 
the above ensures that $m(\cdot)$ is $\nu_{\scriptscriptstyle G(\nu_0, \cdot)}$-a.s.\ the \emph{derivative of the signed measure $\nu_{\scriptscriptstyle F(\nu_0, \cdot)}$ with respect to the measure $\nu_{\scriptscriptstyle G(\nu_0, \cdot)}$}, denoted by $d\nu_{\scriptscriptstyle F(\nu_0, \cdot)}/d\nu_{\scriptscriptstyle G(\nu_0, \cdot)}$, or of  $F(\nu_0, \cdot)$ with respect to $G(\nu_0, \cdot)$, denoted by $d F(\nu_0, \cdot)/d G(\nu_0, \cdot)$.

\subsection{Relations between marginal project metrics.}
\label{s:fpupcli}
This section presents relations between  marginal metrics under threshold policies as functions of the threshold.
The following result gives counterparts of Lemmas \ref{lma:phimpgint} and \ref{lma:pcli3inf} for marginal  metrics $f(x, \cdot)$ and $g(x, \cdot)$. Note that the LS integrals below are well defined, since $g(x, \cdot)$ is a function of bounded variation on $\mathbb{R}$ (see Lemma \ref{pro:gbv}) and $m(\cdot)$ is continuous under PCL-indexability.
The result refers to the finite signed measure $\nu_{g(x, \cdot)}$ determined by $g(x, \cdot)$ through $\nu_{g(x, \cdot)}(z_1, z_2] \triangleq g(x, z_1) - g(x, z_2)$ for finite $z_1 < z_2$.

\begin{lemma}
\label{lma:asspcliii} Under PCL-indexability$,$ 
\begin{itemize}
\item[\textup{(a)}] $m(\cdot)$ is $\nu_{g(x, \cdot)}$-integrable$,$ i.e.$,$ 
the integral $\int |m|  \, d|\nu_{g(x, \cdot)}|$
is finite$;$
\item[\textup{(b)}]
$
f(x, z_2) - f(x, z_1) = \int_{(z_1, z_2]} m(z)  \, g(x, dz),$ for $-\infty \leqslant z_1 < z_2 \leqslant \infty.
$
\end{itemize}
\end{lemma}
\proof (a) This part follows from Lemma \ref{lma:phimpgint} using that, by  (\ref{eq:fxBgxB}), $\nu_{g(x, \cdot)} = \beta \diff_{a=0}^{a=1} \nu_{\scriptscriptstyle G(\kappa^{\scriptscriptstyle a}(x, \cdot), \cdot)}$.

(b)
Using part (a), Lemma \ref{lma:pcli3inf} and (\ref{eq:fxBgxB}), we obtain
\begin{align*}
\int_{(z_1, z_2]} m(z)  \, g(x, dz) & = 
\beta \diff\limits_{a=0}^{a=1}  \int_{(z_1, z_2]} m(z)  \, G(\kappa^a(x, \cdot), dz) \\
& = \beta \diff\limits_{a=0}^{a=1} \{F(\kappa^a(x, \cdot), z_2) - F(\kappa^a(x, \cdot), z_1)\} = f(x, z_2) - f(x, z_1).
\end{align*}
 \endproof

\section{Further relations between MP metrics.}
\label{s:opicfpmpm}
This section  builds on the above to present further relations between MP metrics that will play a key role in the proof of the verification theorem.
The next result shows that the MP index satisfies certain linear Volterra--Stieltjes integral equations.

\begin{lemma}
\label{lma:nusxydec}  Suppose the project is PCL-indexable and let $x \in \mathsf{X}$ and $z \in \overline{\mathbb{R}}.$ Then
\hspace{1in}
\begin{itemize}
\item[{\rm (a)}] 
$m(x, z) - m(z)  = 
\int_{z < b <  x} \frac{g(x, b)}{g(x, z)} \,
m(db),$ \quad if $x > z;$ 
\item[{\rm (b)}]
$m(x, z) - m(z)  = - 
\int_{x \leqslant b \leqslant z} \frac{g(x, b)}{g(x, z)} \,
m(db),$ \quad if $x \leqslant z;$
\end{itemize}
\end{lemma}
\proof
(a)  From Lemmas \ref{pro:cadlagFG}, \ref{lma:asspcliii}, and the dominated convergence theorem, 
we have
\[
f(x, x^{\scriptscriptstyle -}) - f(x, z) =
\lim_{y \nearrow x} \, f(x, y) - f(x, z) =
\lim_{y \nearrow x} \, \int_{(z, y]} m(b) \, g(x, db) =
\int_{(z, x)} m(b) \, g(x, db).
\]
Consider the case $z > -\infty$.
Using  the latter identity, integration by parts, and Lemma \ref{pro:mafgids}(c) gives
\begin{equation}
\label{eq:1fxzphigxz}
\begin{split}
f(x, z) & =
f(x, x^{\scriptscriptstyle -}) - \int_{(z, x)} m(b) \, g(x, db) 
 \\
& = f(x, x^{\scriptscriptstyle -})  - 
[m(x) g(x, x^{\scriptscriptstyle -}) - m(z) g(x, z)] + 
\int_{(z, x)} g(x, b) \, m(db) \\
&=  m(z) g(x, z)  +
\int_{(z, x)} g(x, b) \, m(db).
\end{split}
\end{equation} 
The result now follows by dividing each term by $g(x, z)$.

In the case $z = -\infty$, the result follows from the limiting argument
\begin{align*}
f(x, -\infty) - m(-\infty) g(x, -\infty) & = 
\lim_{z' \to -\infty} f(x, z') - m(z') g(x, z')  \\
& = \lim_{z' \to -\infty} \int_{(z', x)} g(x, b) \, m(db) = \int_{(-\infty, x)} g(x, b) \, m(db),
\end{align*}
where we have used  Lemma \ref{pro:cadlagFG}(b), (PCLI1), (\ref{eq:1fxzphigxz}) and (PCLI2).

(b) Consider the case $z < \infty$. Arguing as in part (a), we  obtain
\begin{equation}
\label{eq:fxzphigxz}
\begin{split}
f(x, z) &  
 = f(x, x^{\scriptscriptstyle -}) + \int_{[x, z]}  m(b) 
   \, g(x, db)  \\
& = f(x, x^{\scriptscriptstyle -}) + m(z) g(x, z) - 
   m(x) g(x, x^{\scriptscriptstyle -}) - 
\int_{[x, z]} g(x, b) \,
m(db) \\
& =  m(z) g(x, z)  - \int_{[x, z]} g(x, b) \, m(db).
\end{split}
\end{equation} 

In the case $z = \infty$, the result follows from the limiting argument
\begin{align*}
f(x, \infty) - m(\infty) g(x, \infty) & = 
\lim_{z' \to \infty} f(x, z') - m(z') g(x, z')  \\
& = -  \lim_{z' \to \infty} \int_{[x, z']} g(x, b) \, m(db) = - \int_{[x, \infty)} g(x, b) \, m(db),
\end{align*}
where we have used  Lemma \ref{pro:cadlagFG}(c), (PCLI1), (\ref{eq:fxzphigxz}) and (PCLI2). This completes the proof.
 \endproof

The next result shows that $m(x, z) - m(z)$ and $m(x) - m(z)$ have the same sign.

\begin{lemma}
\label{lma:maxinf}  Suppose the project is PCL-indexable$.$ Then, for any state $x$ and threshold $z,$
\[
\sgn(m(x, z) - m(z)) = \sgn(m(x) - m(z)).
\]
\end{lemma}
\proof
If $x > z$, we use Lemma \ref{lma:nusxydec}(a) and (PCLI1, PCLI2) to obtain
\begin{align*}
\sgn(m(x, z) - m(z))  & = \sgn
\int_{z < b < x} g(x, b) \,
m(db) =  \sgn(m(x) - m(z)) \geqslant 0.
\end{align*}

If $x \leqslant z$, we use Lemma \ref{lma:nusxydec}(b) and (PCLI1, PCLI2) to obtain
\begin{align*}
\sgn(m(x, z) - m(z))  & = -\sgn  
\int_{x \leqslant b \leqslant z} g(x, b) \,
m(db)  =  \sgn(m(x) - m(z)) \leqslant 0.
\end{align*}
 \endproof

We need a final preliminary result, showing that $m(x, z) - \lambda$ and $m(x) - \lambda$ have the same sign when threshold $z$ is taken from the  set $\mathsf{Z}_\lambda$ defined in (\ref{eq:tslambm}).
 
\begin{lemma}
\label{lma:dpbanix1}
Suppose the project is PCL-indexable$.$ Then$,$ for any $x \in \mathsf{X},$ $\lambda \in \mathbb{R}$ and $z \in \mathsf{Z}_\lambda,$
\begin{equation}
\label{eq:sgnsmxzlmxl}
\sgn(m(x, z) - \lambda) = \sgn(m(x) - \lambda).
\end{equation}
\end{lemma}
\proof Let $x \in \mathsf{X},$ $\lambda \in \mathbb{R}$ and $z \in \mathsf{Z}_\lambda$.
We distinguish three cases. 
In the case $\{y \in \mathsf{X}\colon m(y) = \lambda\} \neq \emptyset$,
we have  $m(z) = \lambda$ (see (\ref{eq:tslambm})), and hence (\ref{eq:sgnsmxzlmxl}) holds by Lemma \ref{lma:maxinf}.

If $\lambda < m(y)$ for every $y$, then (see (\ref{eq:tslambm})) $z < \ell$ if $\ell > -\infty$, and $z = -\infty$ if $\ell = -\infty$, so  $m(z) = m(\ell) \geqslant \lambda$, with strict inequality if $\ell > -\infty$ (as then $\ell \in \mathsf{X}$).
We can now write, using Lemma \ref{lma:nusxydec}(a),
\begin{equation}
\label{eq:mxzml1}
m(x, z) - \lambda = m(x, z) - m(z) + m(z) - \lambda  = \int_{(z, x)} \frac{g(x, b)}{g(x, z)} \, m(db) + m(\ell) - \lambda > 0,
\end{equation}
where the inequality follows from (PCLI1) if $m(\ell) > \lambda$. If  $m(\ell) = \lambda$ then 
$\ell = -\infty$, and by the assumptions the integral in the right-hand side of (\ref{eq:mxzml1}) is positive, yielding the inequality and  (\ref{eq:sgnsmxzlmxl}).

If $\lambda > m(y)$ for every $y$, then (see (\ref{eq:tslambm})) $z > u$ if $u < \infty$, and $z = \infty$ if $u = \infty$, so $m(z) = m(u) \leqslant \lambda$, with strict inequality if $u$ is finite (because then $u \in \mathsf{X}$).
We can now write, using Lemma \ref{lma:nusxydec}(b),
\begin{equation}
\label{eq:mxzml2}
m(x, z) - \lambda = m(x, z) - m(z) + m(z) - \lambda  = -\int_{b \geqslant x} \frac{g(x, b)}{g(x, z)} \, m(db) + m(u) - \lambda < 0,
\end{equation}
where the inequality follows from (PCLI1) if $m(u) < \lambda$. If  $m(u) = \lambda$ then 
$u = \infty$, and by the assumptions the integral in the right-hand side of (\ref{eq:mxzml2}) is positive, which yields (\ref{eq:sgnsmxzlmxl}).
 \endproof

\section{Proving Theorem \ref{the:pcliii}.}
\label{s:ptpcliii}
We are now ready to prove Theorem \ref{the:pcliii}.

\proof{Proof of Theorem $\ref{the:pcliii}.$}
Let $\zeta(\cdot)$ be a generalized inverse of the MP index $m(\cdot)$. 
For any $x \in \mathsf{X}$ and $\lambda \in \mathbb{R}$, we have $\zeta(\cdot) \in \mathsf{Z}_\lambda$ and, using (PCLI1) and Lemma \ref{eq:sgnsmxzlmxl}, 
\begin{equation}
\label{eq:sgnproofth}
\sgn \Delta_{a=0}^{a=1} V_\lambda(x, z) = \sgn (f(x, \zeta(\lambda)) - \lambda g(x, \zeta(\lambda))) 
                                        = \sgn (m(x, \zeta(\lambda)) - \lambda) 
                                        = \sgn (m(x) - \lambda).
\end{equation}
The outermost identity in (\ref{eq:sgnproofth}) and Lemma \ref{lma:geninvprop}(a) imply that 
 $\zeta(\cdot)$ satisfies (\ref{eq:thrlsigndef}), whence the project is thresholdable with (see Remark \ref{re:thrlsigndef}(i))
optimal value $V_\lambda^*(x) = V_\lambda(x, \zeta(\lambda))$. We thus obtain 
\begin{equation}
\label{eq:finalsgnproofth}
\sgn \Delta_{a=0}^{a=1} V_\lambda^*(x) = \sgn (m(x) - \lambda),
\end{equation}
which shows (cf.\ Definition \ref{def:indx}) that the project is indexable with Whittle index $m(\cdot)$.
 \endproof

\section{How to establish condition (PCLI3).}
\label{s:ccpcli3} 
This section presents conditions under which PCL-indexability condition
 (PCLI3) can be established by simpler means than direct verification. 

\subsection{Piecewise differentiable $F(x, \cdot)$ and $G(x, \cdot)$.}
\label{s:tcopwdiffG}
We start with the case that metrics $F(x, \cdot)$ and $G(x, \cdot)$ are piecewise differentiable.
Denote by  $F_2'(x, z)$ and $G_2'(x, z)$ their partial derivatives with respect to $z$.
Below, $\mathring{S}$ denotes the interior of a set $S$.

\begin{assumption}
\label{ass:pwdiffFG} For each  state $x$ there is a countable partition $\{S_i(x)\colon i \in N(x)\}$ of  $\mathbb{R}$ with each $S_i(x)$ being either a left-semiclosed interval or a singleton$,$ such that$:$
\begin{itemize}
\item[\textup{(i)}]
$F(x, \cdot) = \sum_{i \in N(x)} F_i(x, \cdot) 1_{S_i(x)}(\cdot)$ and $G(x, \cdot) = \sum_{i \in N(x)} G_i(x, \cdot) 1_{S_i(x)}(\cdot),$ with $F_i(x, \cdot)$ and $G_i(x, \cdot)$ being both differentiable  on  $\mathring{S}_i(x).$ 
\item[\textup{(ii)}]  $F_2'(x, z) = m(z) G_2'(x, z)$ for $z \in \mathring{S}_i(x)$ and $i \in N(x).$ 
\end{itemize}
\end{assumption}

We need two preliminary results.
Denote by $\mathrm{TV}_{\scriptscriptstyle F(x, \cdot)}(I)$ and $\mathrm{TV}_{\scriptscriptstyle G(x, \cdot)}(I)$ the \emph{total variation} of $F(x, \cdot)$ and $G(x, \cdot)$ on an interval $I$, and  by $D(x) \triangleq \mathbb{R} \setminus \cup_{i \in N(x)} \mathring{S}_i(x)$ the set of interval endpoints for  partition $\{S_i(x)\colon i \in N(x)\}$. 
Recall that $\mathbb{V}(I)$ denotes the functions of bounded variation on $I$.

\begin{lemma} 
\label{lma:pwdiffFG}
Let \textup{(PCLI1, PCLI2)} and Assumption $\ref{ass:pwdiffFG}$ hold$.$
Let $x$ be a state and $I$ a finite open interval  where both
$F(x, \cdot)$ and $G(x, \cdot)$ are differentiable$.$ Then $F(x, \cdot) \in \mathbb{V}(I)$  and
\[
\int_I F(x, dz) = \int_I m(z) \, G(x, dz).
\]
\end{lemma}
\proof
Since $G(x, \cdot)$ is bounded nonincreasing by Lemmas \ref{pro:cadlagFG}(a) and \ref{pro:Gnondecr}(a), and $m(\cdot)$ is continuous, by standard results $\int_I |m|(z) \, G(x, dz)$ is finite.
Furthermore,  using Assumption \ref{ass:pwdiffFG}(ii),
\begin{align*}
\int_I |m|(z) \, G(x, dz) & = \int_I |m|(z) G_2'(x, z) \, dz = -\int_I |F_2'|(x, z) \, dz = -\mathrm{TV}_{\scriptscriptstyle F(x, \cdot)}(I)
\end{align*}
(cf.\ \citet[theorem 6.1.7]{cartvanBrunt00}). Hence, $F(x, \cdot) \in \mathbb{V}(I)$ and \begin{align*}
\int_I m(z) \, G(x, dz) & = \int_I m(z) G_2'(x, z) \, dz = \int_I F_2'(x, z) \, dz = \int_I F(x, dz).
\end{align*}
\qquad
 \endproof

In the next result, 
$\mathbb{V}_{\textup{loc}}(\mathbb{R})$ denotes the class of functions of \emph{locally bounded variation} on $\mathbb{R}$, i.e., having bounded variation on every finite interval.

\begin{lemma} 
\label{lma:pwdiffFG2}
Under \textup{(PCLI1, PCLI2)} and Assumption $\ref{ass:pwdiffFG},$ $F(x, \cdot) \in \mathbb{V}_{\textup{loc}}(\mathbb{R})$.
\end{lemma}
\proof
Let $I = (z_1, z_2)$ with $-\infty < z_1  < z_2 < \infty$ and $K(I) \triangleq \max \{|m|(z)\colon z \in [z_1, z_2]\}$, which is finite by  continuity of $m(\cdot)$.
The total variation $\mathrm{TV}_{\scriptscriptstyle F(x, \cdot)}(I)$ of $F(x, \cdot)$ on $I$ is finite, since
\begin{align*}
\mathrm{TV}_{\scriptscriptstyle F(x, \cdot)}(I) & = \sum_{z \in D(x) \cap I} |\Delta_{2} F(x, z)| + \sum_{i\colon \mathring{S}_i(x) \cap I \neq \emptyset} \int_{\mathring{S}_i(x) \cap I} |F_2'|(x, z) \, dz   \\
& \quad =
\sum_{z \in D(x) \cap I} |m|(z) |\Delta_{2} G(x, z)| + \sum_{i\colon \mathring{S}_i(x) \cap I \neq \emptyset} \int_{\mathring{S}_i(x) \cap I} |m|(z) |G_2'|(x, z) \, dz \\ 
& \quad \leqslant 
K(I) \bigg\{\sum_{z \in D(x) \cap I} |\Delta_{2} G(x, z)|  + \int_{i\colon \mathring{S}_i(x) \cap I \neq \emptyset} |G_2'|(x, z) \, dz\bigg\}  \\
& \quad = K(I) \, \mathrm{TV}_{\scriptscriptstyle G(x, \cdot)}(I) \leqslant K(I) \, \mathrm{TV}_{\scriptscriptstyle G(x, \cdot)}(\mathbb{R}) < \infty,
\end{align*}
where we have used Assumption \ref{ass:pwdiffFG} and Lemmas \ref{pro:mafgids}(a),  \ref{lma:pwdiffFG} and \ref{pro:gbv}.
\qquad
 \endproof

We can now give the main result of this section.

\begin{proposition} 
\label{pro:pwdiffFG}
Under \textup{(PCLI1, PCLI2)} and Assumption $\ref{ass:pwdiffFG},$ \textup{(PCLI3)} holds$.$
\end{proposition}
\proof
Let $-\infty < z_1 < z_2 < \infty$. 
We have
\begin{align*}
& F(x, z_2) - F(x, z_1) = \int_{(z_1, z_2]} F(x, dz) = \Delta_{2} F(x, z_2) + \int_{(z_1, z_2)} F(x, dz) \\
& \quad = \sum_{z \in D(x) \cap (z_1, z_2]} \Delta_{2} F(x, z) + \sum_{i\colon \mathring{S}_i(x) \cap (z_1, z_2) \neq \emptyset} \int_{\mathring{S}_i(x) \cap (z_1, z_2)} F_2'(x, z) \, dz   \\
& \quad =
\sum_{z \in D(x) \cap (z_1, z_2]} m(z) \, \Delta_{2} G(x, z) + \sum_{i\colon \mathring{S}_i(x) \cap (z_1, z_2) \neq \emptyset} \int_{\mathring{S}_i(x) \cap (z_1, z_2)} m(z) \, G_2'(x, z) \, dz \\ 
& \quad = 
\int_{(z_1, z_2]} m(z) \, G(x, dz),
\end{align*}
where we have used Lemma \ref{lma:pwdiffFG2}, Assumption \ref{ass:pwdiffFG} and Lemmas \ref{pro:mafgids}(a),  \ref{lma:pwdiffFG} and \ref{lma:eaafpcli3}.
\qquad
 \endproof

\subsection{Piecewise constant $G(x, \cdot)$.}
\label{s:tcopwconstG}
We next consider the case that the resource metric $G(x, z)$ is piecewise constant in the threshold variable $z$, with a finite or denumerable number of pieces.

\begin{assumption}
\label{ass:pwcFG}  For each state $x,$ 
$G(x, \cdot) = \sum_{i \in N(x)} G_i(x) 1_{S_i(x)}(\cdot),$ where $\{S_i(x)\colon i \in N(x)\}$ is a countable partition of $\mathbb{R}$ with each $S_i(x)$ being either a left-semiclosed interval or a singleton. 
\end{assumption}

Note that Assumption \ref{ass:pwcFG} does not refer to the reward metric, so it is not evident that this case reduces to that  in Section \ref{s:tcopwdiffG}. Yet, we will show
that Assumption \ref{ass:pwcFG} implies Assumption  \ref{ass:pwdiffFG}.

The following result can be used to verify Assumption \ref{ass:pwcFG}  through analysis of state trajectories $\{X_t\}_{t=0}^\infty$  under threshold policies. Parts (a, c) characterize when $G(x, \cdot)$ is constant on an interval. Parts (b, d) ensure that  $F(x, \cdot)$ is  constant over the same intervals that $G(x, \cdot)$.
For a Borel set $S$, denote by
$T_S \triangleq \min\{t \geqslant 0\colon X_t \in S\}$ the first hitting time to $S$, with $T_S \triangleq \infty$ if  $S$ is never hit. 

\begin{lemma}
\label{pro:GcImplFc}
Let  \textup{(PCLI1)} hold$.$ For any state $x$ and finite $z_1 < z_2,$
\begin{itemize}
\item[\textup{(a)}] $G(x, \cdot)$ is constant on $[z_1, z_2]$ iff $\Prob_{x}^{z_1}\{T_{(z_1, z_2]} < \infty\} = 0;$ 
\item[\textup{(b)}]
 if $G(x, \cdot)$ is constant on  $[z_1, z_2],$ then so is $F(x, \cdot);$
\item[\textup{(c)}] $G(x, \cdot)$ is constant on $[z_1, z_2)$ iff $\Prob_{x}^{z_1}\{T_{(z_1, z_2)} < \infty\} = 0;$
\item[\textup{(d)}]
 if $G(x, \cdot)$ is constant on  $[z_1, z_2),$ then so is $F(x, \cdot).$
\end{itemize}
\end{lemma}
\proof
(a)
Since $G(x, \cdot)$ is nonincreasing on $\mathbb{R}$ by Lemma \ref{pro:Gnondecr}(a),  it is constant on $[z_1, z_2]$ iff $G(x, z_1) =
  G(x, z_2)$, which occurs,
by Lemma \ref{lma:FGxzdfgxz}(b), iff
\begin{equation}
\label{eq:Gxz1z2}
\int_{(z_1, z_2]} g(y, z_2) \,
  \tilde{\mu}_{x}^{z_1}(dy) = 0.
  \end{equation}
Now, under (PCLI1), (\ref{eq:Gxz1z2}) holds iff $\tilde{\mu}_{x}^{z_1}(z_1, z_2] = 0$, which occurs iff $\Prob_{x}^{z_1}\{T_{(z_1, z_2]} < \infty\} = 0.$

 (b) 
If $G(x, \cdot)$ is constant on $(z_1, z_2]$ then $\tilde{\mu}_{x}^{z_1}(z_1, z_2] = 0$, which, along with  Lemma \ref{lma:FGxzdfgxz}(a), yields 
\[
F(x, z_1) -
  F(x, z) = \int_{(z_1, z]} f(y, z) \,
  \tilde{\mu}_{x}^{z_1}(dy)  = 0, \quad z \in (z_1, z_2].
\]
Hence, $F(x, \cdot)$ is also constant on  $[z_1, z_2]$. 

(c) Let $\{z^n\}_{n=0}^\infty$ be increasing with $z_1 < z^n < z_2$ and 
$\lim_{n \to \infty} z^n = z_2$. Then 
$G(x, \cdot)$ is constant on $[z_1, z_2)$ iff $G(x, \cdot)$ is constant on $[z_1, z^n]$ for each $n$, which by part (a) occurs  iff $\tilde{\mu}_{x}^{z_1}(z_1, z^n] = 0$ for each $n$. By continuity from below,  such will be the case iff $\tilde{\mu}_{x}^{z_1}(z_1, z_2) = 0$, i.e.,  iff $\Prob_{x}^{z_1}\{T_{(z_1, z_2)} < \infty\} = 0$. 

(d) The proof of this part follows along the lines of that for part (b).
 \qquad
 \endproof

Note that, e.g., part (c) says that $G(x, \cdot)$ is constant on $[z_1, z_2)$ iff the project state never hits $(z_1, z_2)$ under the $z_1$-policy starting from $x$.
Parts (b, d) ensure that  $F(x, \cdot)$ is of the form $F(x, \cdot) = \sum_{i \in N(x)} F_i(x) 1_{S_i(x)}(\cdot)$.

\begin{lemma}
\label{lma:as1ias2}
Assumption $\ref{ass:pwcFG}$ implies Assumption  $\ref{ass:pwdiffFG}.$
\end{lemma}
\proof
The result follows immediately from Lemma \ref{pro:GcImplFc}(d).
 \endproof

Verification of Assumption \ref{ass:pwcFG} for a given model may be a nontrivial task, as it requires knowledge of $G(x, \cdot)$.
Instead, we will draw on Lemma \ref{pro:GcImplFc}
to show that it suffices to verify 
the following assumption, which only involves analyzing project state trajectories under threshold policies.

\begin{assumption}
\label{ass:pwctraj}  For each $x$ there is a countable set  $D(x)$ giving the endpoints of a partition of $\mathbb{R}$ into left-semiclosed   intervals or singletons$,$ such that $\mathsf{P}_x^z\{\{X_t
\}_{t=0}^\infty \subseteq D(x)\} = 1$ for any $z.$  
\end{assumption}

Note that Assumption \ref{ass:pwctraj} ensures that, starting from $X_0 = x$, the project state trajectory $\{X_t\}_{t = 0}^\infty$ lies in $D(x)$ under \emph{any} threshold policy.

\begin{lemma}
\label{lma:pwctraj}
Assumption $\ref{ass:pwctraj}$ implies Assumption  $\ref{ass:pwcFG}.$
\end{lemma}
\proof
The result follows immediately from Lemma \ref{pro:GcImplFc}(c), using the points in $D(x)$ to define the endpoints of the partition $\{S_i(x)\colon i \in N(x)\}$.
 \endproof

We thus obtain the following result.

\begin{proposition}
\label{pro:1cpcli3} Under \textup{(PCLI1, PCLI2)} and Assumption $\ref{ass:pwcFG}$ or $\ref{ass:pwctraj},$ \textup{(PCLI3)} holds$.$
\end{proposition}
\proof  
The result follows directly from Lemmas \ref{lma:as1ias2} and \ref{lma:pwctraj} and Proposition \ref{pro:pwdiffFG}.
\qquad
 \endproof

\section{Metrics and index computation.}
\label{s:pmaic}
This section gives recursions for approximately computing the project metrics and index, in cases they  cannot be evaluated in closed form.

For a state $x$, threshold $z$ and $k = 0, 1, \ldots,$
consider the \emph{$k$-horizon metrics}
\begin{equation}
\label{eq:fhFGboldx0pis}
F_k(x, z) \triangleq \Ex_{x}^z\left[\sum_{t = 0}^k 
\beta^t r(X_{t}, A_{t})\right] \quad \textup{and} \quad G_k(x, z) \triangleq \Ex_{x}^z\Bigg[\sum_{t = 0}^k 
\beta^t c(X_{t}, A_{t})\Bigg].
\end{equation}
The function sequences $\{F_{k}(\cdot, z)\}_{k=0}^\infty$ and $\{G_{k}(\cdot, z)\}_{k=0}^\infty$ are determined by the following \emph{value iteration} recursions: $F_0(x, z) \triangleq r(x, 1_{\{x > z\}})$,  $G_0(x, z) \triangleq c(x, 1_{\{x > z\}})$ and, for $k = 0, 1, \ldots$, 
\begin{equation}
\label{eq:viFk}
F_{k+1}(x, z) \triangleq 
r(x, 1_{\{x > z\}}) + \beta \int F_{k}(y, z) \, \kappa^{1_{\{x > z\}}}(x, dy), 
\end{equation}
\begin{equation}
\label{eq:viGk}
G_{k+1}(x, z) \triangleq 
c(x, 1_{\{x > z\}}) + \beta \int G_{k}(y, z) \, \kappa^{1_{\{x > z\}}}(x, dy).
\end{equation}

Consider further  the \emph{$k$-horizon marginal metrics} $f_k(x, z)$ and $g_k(x, z)$, given by $f_0(x, z)~=~\Delta_{a=0}^{a=1} r(x, a)$, $g_0(x, z) = \Delta_{a=0}^{a=1} c(x, a)$ and, for $k \geqslant 1$, 
\begin{equation}
\label{eq:fkxz}
f_{k}(x, z) = \Delta_{a=0}^{a=1} \{r(x, a) + \beta \int F_{k-1}(y, z) \, \kappa^{a}(x, dy)\},
\end{equation}
\begin{equation}
\label{eq:gkxz}
g_{k}(x, z) = \Delta_{a=0}^{a=1} \{c(x, a) + \beta \int G_{k-1}(y, z) \, \kappa^{a}(x, dy)\}.
\end{equation}

Finally, consider the \emph{$k$-horizon MP metric} $m_k(x, z) \triangleq f_{k}(x, z)/g_{k}(x, z)$ and the \emph{$k$-horizon MP index} $m_k(x) \triangleq m_k(x, x)$, which are defined when $g_{k}(x, z) \neq 0$ and $g_{k}(x, x) \neq 0$, respectively.

The next result shows, in parts (a, c), that $F_{k}(\cdot, z)$, $G_{k}(\cdot, z)$, $f_{k}(\cdot, z)$ and $g_{k}(\cdot, z)$ converge linearly with rate $\gamma$ in the $w$-norm $\|\cdot\|_w$ (see (\ref{eq:vxwnorm}) and Remark \ref{re:assboundfs}(ii)). 
 Parts (b, d) bound their  $w$-norms.
Note that $M$, $\gamma$ and $w$ are as in Assumption \ref{ass:first}(ii), and $M_\gamma$ as in (\ref{eq:FGMwnon}).

\begin{lemma}
\label{lma:cfhpm} For any threshold $z$ and integer $k \geqslant 0,$
\begin{itemize}
\item[\textup{(a)}] $\|F_{k}(\cdot, z) - F(\cdot, z)\|_w  \leqslant M_\gamma  \gamma^k$ and $\|G_{k}(\cdot, z) - G(\cdot, z)\|_w  \leqslant M_\gamma  \gamma^k;$
\item[\textup{(b)}] $\|F_{k}(\cdot, z)\|_w  \leqslant M_\gamma$ and $\|G_{k}(\cdot, z)\|_w  \leqslant M_\gamma;$
\item[\textup{(c)}] $\|f_{k}(\cdot, z) - f(\cdot, z)\|_w  \leqslant 2 M_\gamma  \gamma^k$ and $\|g_{k}(\cdot, z) - g(\cdot, z)\|_w  \leqslant 2 M_\gamma  \gamma^k;$
\item[\textup{(d)}] $\|f_{k}(\cdot, z)\|_w  \leqslant 2 M_\gamma$ and $\|g_{k}(\cdot, z)\|_w  \leqslant 2 M_\gamma;$
\end{itemize}
\end{lemma}
\proof We only consider reward metrics, as the results for resource metrics follow similarly.

(a) (\ref{eq:viFk}) is a contraction with modulus $\gamma$  on  $\mathbb{B}_w(\mathsf{X})$ (cf.\  \citet[remark 8.3.10]{herlerLass99}), so 
 (cf.\ Remark \ref{re:assboundfs}(iv)) Banach's fixed point theorem yields the result.

(b) We first show by induction on $k$ that
$|F_k|(x, z) \leqslant C_k w(x)$, where $C_0 \triangleq 0$ and 
$C_k \triangleq M (1 + \cdots + \gamma^{k-1})$ for $k \geqslant 1$, so $C_{k+1} = M  + \gamma C_k$.
For $k = 0$, such a result trivially holds.  Assume now $|F_k|(x, z) \leqslant C_k w(x)$ for some $k \geqslant 0$. Then, from Assumption \ref{ass:first} and (\ref{eq:viFk}) we obtain
\begin{align*}
|F_{k+1}|(x, z)
& \leqslant |r|(x, 1_{\{x > z\}}) + \beta \int |F_{k}|(y, z) \, \kappa^{1_{\{x > z\}}}(x, dy)  \leqslant  M w(x) + \beta \int C_k w(y) \, \kappa^{1_{\{x > z\}}}(x, dy) \\
& \leqslant  M w(x) + \gamma C_k w(x) = C_{k+1} w(x),
\end{align*}
which proves the induction step.
Since $C_k < M_\gamma$, it follows that $|F_k|(x, z) \leqslant M_\gamma w(x).$

(c) Using (\ref{eq:mrm1}), (\ref{eq:fkxz}), part (a) and Assumption \ref{ass:first}(ii.b), we obtain
\begin{align*}
 |f_{k}(x, z) - f(x, z)| & \leqslant \beta \int |F_{k-1}(y, z) - F(y, z)| \, \kappa^{1}(x, dy) + \beta \int |F_{k-1}(y, z) - F(y, z)| \, \kappa^{0}(x, dy) \\
 & \leqslant \beta M_\gamma \gamma^{k-1} \int   w(y)  \, \kappa^{1}(x, dy) + \beta M_\gamma \gamma^{k-1} \int   w(y)  \, \kappa^{0}(x, dy) \leqslant 2 M_\gamma  \gamma^{k} w(x).
 \end{align*} 
 
 (d) This part follows along the same lines as part (b).
\qquad
 \endproof

Letting $\underline{g}(x, z) \triangleq \inf_k g_k(x, z)$,  
consider the following strong version of condition (PCLI1):
\begin{equation}
\label{eq:fhcpcli1}
\underline{g}(x, z) > 0 \quad \textup{for any state } x \textup{ and threshold } z.\end{equation}
Since $g_k(x, z)$ converges to $g(x, z)$ (see Lemma \ref{lma:cfhpm}(c)),  condition (\ref{eq:fhcpcli1}) implies (PCLI1).

The following result shows that the finite-horizon MP metric $m_k(x, z)$ and MP index $m_{k}(x)$ converge at least linearly with rate $\gamma$ to  $m(x, z)$ and  $m(x)$, respectively.

\begin{proposition}
\label{pro:cfhMPI} Let \textup{(\ref{eq:fhcpcli1})} hold$.$ Then$,$ for any state $x,$ threshold $z$ and integer $k \geqslant 0,$
\begin{itemize}
\item[\textup{(a)}] $|m_k(x, z) - m(x, z)|  \leqslant 2 M_\gamma   \gamma^k w(x) (1+|m|(x, z)) /\underline{g}(x, z);$
\item[\textup{(b)}] $|m_{k}(x) - m(x)|  \leqslant 2 M_\gamma   \gamma^k w(x) (1+|m|(x)) /\underline{g}(x, x).$
\end{itemize}
\end{proposition}
\proof 
(a) From
\[
m_k(x, z) - m(x, z) = \frac{f_{k}(x, z) - f(x, z)}{g_{k}(x, z)} - \frac{g_{k}(x, z) - g(x, z)}{g_{k}(x, z)} m(x, z),
\]
Lemma \ref{lma:cfhpm}(c) and $g_{k}(x, z) \geqslant \underline{g}(x, z) > 0$, we obtain the stated inequality.

(b) This part follows by setting $z = x$ in part (a).
\qquad
 \endproof

\section{Examples.}
\label{s:ex}
This section illustrates application of the proposed approach to two  models.

\subsection{Optimal web crawling.}
\label{s:owcec}
This is the model analyzed in 
\citet{avraBorkar18} under the average-reward criterion, which we analyze below under the discounted criterion. The project has reward and resource consumption 
$r(x, a) \triangleq x a$ and $c(x, a) \triangleq C a$ with $C > 0$, and the following dynamics. When $X_t = x$ and $A_t = 0$, $X_{t+1} = \ell + \alpha x$, for certain parameters $0 < \alpha < 1$ and $\ell =  (1-\alpha) b$ with $b > 0$; if $A_t = 1$, $X_{t+1} = \ell$. 
The state space is $\mathsf{X} \triangleq [\ell, u]$, with $u \triangleq \ell/(1-\alpha)$. 

The evaluation equations for the reward and resource metrics under the $z$-threshold policy are
\[
F(x, z) = 
\begin{cases} 
\beta F(\ell + \alpha x, z) & \textup{ if } x \leqslant z \\
x + \beta F(\ell, z) & \textup{ if } x > z
\end{cases}
\qquad 
\textup{and}
\qquad 
G(x, z) = 
\begin{cases} 
\beta G(\ell + \alpha x, z) & \textup{ if } x \leqslant z \\
C + \beta G(\ell, z) & \textup{ if } x > z.
\end{cases}
\]

To solve such equations in closed form,  
let $h(x) \triangleq \ell + \alpha x$, and recursively define $h_0(x) \triangleq x$ and 
$h_t(x) \triangleq h(h_{t-1}(x))$ for $t \geqslant 1$, so 
 $h_t(x) = u - (u - x) \alpha^t$.
 Let $\tau(x, z) \triangleq \min \{t \geqslant 0\colon h_t(x) > z\}$, with $\tau(x, z) \triangleq \infty$ if $h_t(x) \leqslant z$ for all $t$. Through elementary arguments we obtain

\[
F(x, z) = 
\begin{cases} 
\beta^{\tau(x, z)} (h_{\tau(x, z)}(x) + \beta F(\ell, z)) & \textup{ if } x \leqslant z \\
x + \beta F(\ell, z) & \textup{ if } x > z
\end{cases}
\qquad 
\textup{with}
\qquad 
F(\ell, z) = \frac{\beta^{\tau(\ell, z)} h_{\tau(\ell, z)}(\ell)}{1-\beta^{\tau(\ell, z)+1}},
\]
\[
G(x, z) = 
\begin{cases} 
\beta^{\tau(x, z)} (C + \beta G(\ell, z)) & \textup{ if } x \leqslant z \\
C + \beta G(\ell, z) & \textup{ if } x > z
\end{cases}
\qquad 
\textup{with}
\qquad 
G(\ell, z) = \frac{\beta^{\tau(\ell, z)} C}{1-\beta^{\tau(\ell, z)+1}}.
\]
Note that, since the $u$-threshold policy is never active, we have
$F(x, u) \equiv 0$ and $G(x, u) \equiv 0$.

Further, the marginal metrics are given by $f(x, u) = x$, $g(x, u) = C > 0$, and
\[
f(x, z) = x + \beta F(\ell, z) - \beta F(h(x), z) =
x + \beta F(\ell, z) - \beta^{\tau(x, z)} (h_{\tau(x, z)}(x) + \beta F(\ell, z)),
\]
\begin{align*}
g(x, z) & = C + \beta G(\ell, z) - \beta G(h(x), z) =
C + \beta G(\ell, z) - \beta^{\tau(x, z)} (C + \beta G(\ell, z)) \\
& =
(1-\beta^{\tau(x, z)}) (C + \beta G(\ell, z)) > 0.
\end{align*}

As for the MP index $m(\cdot)$, since $\tau(x, x) = 1$ for $\ell \leqslant x < u$, we have, for such $x$ and $t = 1, 2, \ldots$, 
\[
m(x) \triangleq 
\frac{f(x, x)}{g(x, x)} = \frac{x - \beta h(x) + \beta (1-\beta) F(\ell, x)}{(1-\beta) (C + \beta G(\ell, x))} =
\frac{x - \beta h(x) + \beta (1-\beta) F_t(\ell)}{(1-\beta) (C + \beta G_t(\ell))}, \;
h_{t-1}(\ell) \leqslant x < h_{t}(\ell),
\]
with $F_t(\ell) \triangleq \beta^{t} h_{t}(\ell)/(1-\beta^{t+1})$ and 
 $G_t(\ell) \triangleq \beta^{t} C/(1-\beta^{t+1})$.
Further, 
$
m(u) \triangleq 
f(u, u)/g(u, u) = u/C.
$

\begin{proposition}
\label{pro:pclowc}
The above optimal web crawling project model is PCL-indexable.
\end{proposition}
\proof  
It is shown above that $g(x, z) > 0$, hence condition (PCLI1) holds.
The MP index $m(\cdot)$ is  piecewise affine right-continuous, increasing within each piece, since 
\[
\frac{d}{dx} m(x) = \frac{1 - \alpha \beta}{(1-\beta) (C + \beta G_t(\ell))} > 0, \quad
h_{t-1}(\ell) < x < h_{t}(\ell), \quad t = 1, 2, \ldots
\]

As for continuity of $m(\cdot)$, we need to check that, for $t = 1, 2, \ldots$, 
$m(h_t(\ell)^{\scriptscriptstyle -}) = m(h_t(\ell))$, i.e.,  
\[
\frac{h_t(\ell) - \beta h_{t+1}(\ell) + \beta (1-\beta) F_t(\ell)}{(1-\beta) (C + \beta G_t(\ell))}
=
\frac{h_t(\ell) - \beta h_{t+1}(\ell) + \beta (1-\beta) F_{t+1}(\ell)}{(1-\beta) (C + \beta G_{t+1}(\ell))},
\]
which is straightforward to verify algebraically.
We further need to check that $m(h_t(\ell)) \to m(u)$ as $t \to \infty$, which also follows immediately. 
This establishes condition (PCLI2). 

For checking condition (PCLI3), we apply the results in Section \ref{s:tcopwconstG}, since 
the resource metric $G(x, z)$ is piecewise constant in the threshold $z$ for fixed $x$, as it is immediate that, for $z < x$, 
\[
G(x, z) =
\begin{cases}
\displaystyle \frac{\beta^{t+1}}{1-\beta^{t+1}} C, & \quad x > z, h_{t-1}(\ell) \leqslant z < h_t(\ell), \, t \geqslant 1 \\
\displaystyle \frac{\beta^{s+t+1}}{1-\beta^{t+1}} C, & \quad x \leqslant z, h_{s-1}(x) \leqslant z < h_s(x), h_{t-1}(\ell) \leqslant z < h_t(\ell), \, 1 \leqslant s \leqslant t.
\end{cases}
\]

Thus, Assumption \ref{ass:pwcFG} holds, and hence Proposition \ref{pro:1cpcli3} implies that (PCLI3) holds.
\qquad
 \endproof

Hence, Theorem \ref{the:pcliii} gives that the model is threshold-indexable with index $m(\cdot)$. 
Note further that the index $m(\cdot)$ converges as $\beta \to 1$ to the limit
\[
\overline{m}(x) = \frac{(t+1)(x - h(x)) + h_t(\ell)}{C}, \quad 
h_{t-1}(\ell) \leqslant x < h_{t}(\ell), \, t \geqslant 1,
\]
which is precisely the average-reward Whittle index derived in \citet{avraBorkar18}.

\subsection{Optimal dynamic transmission over a noisy channel.}
\label{s:ocam}
This is the model in \citet{liuZhao10}, which motivated the proposal of a previous incomplete version of the PCL-indexability conditions for real-state projects in \citet{nmngi08}.
This section completes the analysis in \citet{nmngi08}, and further proves satisfaction of (PCLI3), which was not considered there. 

A user dynamically attempts to send packets over a noisy channel.
In each period $t$, the channel can be in either state $1$ (good), in which an  transmission attempt succeeds, or $0$, in which it fails. State transitions are Markovian, with probabilities $p$ and $q$ for  transitions $1 \to 0$ and $0 \to 1$, satisfying $\rho \triangleq 1 - p - q > 0$. The user cannot observe the channel state, and bases decisions on the \emph{belief state} $X_t \in \mathsf{X} \triangleq [0, 1]$,  the posterior probability that the channel state is $1$.
If $X_t = x$ and action $A_t = 1$ is taken (attempt to transmit) then $X_{t+1}$ takes the values  $q + \rho$ and $q$ with probabilities $x$ and $1-x$.
If action $A_t = 0$ is taken, $X_{t+1} = q + \rho x$.
The reward and resource functions are $r(x, a) \triangleq a x$ and $c(x, a) \triangleq a$.
The cost per transmission attempt is $\lambda$.

The $\lambda$-price problem $P_\lambda$ in (\ref{eq:lpricepfg1}) is to find a transmission policy  maximizing the expected total $\beta$-discounted value.
To apply the present framework we need to evaluate the project metrics. Under the $z$-policy, the reward and resource metrics are determined by the functional equations
\begin{align*}
F(x, z) & = 
\beta F(q + \rho x, z) 1_{\{x \leqslant z\}} + (x + \beta x F(q + \rho, z) + \beta (1-x) F(q, z)) 1_{\{x > z\}}, \\
G(x, z) & = 
\beta G(q + \rho x, z) 1_{\{x \leqslant z\}} +
(1 + \beta x G(q + \rho, z) + \beta (1-x) G(q, z)) 1_{\{x > z\}},
\end{align*}
and the marginal metrics are given by
\begin{align*}
f(x, z) & = x - \beta (F(q + \rho x, z)-x F(q + \rho, z) - (1-x) F(q, z)), \\
g(x, z) & = 
1 - \beta (G(q + \rho x, z)-x G(q + \rho, z) - (1-x) G(q, z)).
\end{align*}
Provided $g(x, x) \neq 0$, the MP index is 
\[
m(x) = \frac{x - \beta (F(q + \rho x, x)-x F(q + \rho, x) - (1-x) F(q, x))}{1 - \beta (G(q + \rho x, x)-x G(q + \rho, x) - (1-x) G(q, x))}.
\]

The following analysis 
uses the map $h(x) \triangleq q + \rho x$ and \emph{iterates} $h_t(x)$, with $h_0(x) = x$ and 
$h_t(x) = h(h_{t-1}(x))$ for $t \geqslant 1$, so  $h_t(x) = h_{\infty} - (h_{\infty} - x) \rho^t \to h_\infty \triangleq q/(1-\rho)$ as $t \to \infty$.
It further uses the \emph{backward iterates} $h_{-t}(z) = h_{\infty} - (h_{\infty} - z) \rho^{-t}$ for $t \geqslant 1$.
Four cases need be considered.

\paragraph{Case I: $z < q$.} In this case  the belief state is above threshold at times $t \geqslant 1$. The metrics are
\begin{align*}
F(x, z) & = 
\frac{\beta  (q+(1-\beta ) \rho  x) + (1-\beta) (1-\beta  \rho) x  1_{\{x > z\}}}{(1-\beta) (1-\beta \rho)}, & f(x, z) & = x \\
 G(x, z) & = \frac{\beta + (1-\beta) 1_{\{x > z\}}}{1-\beta}, & g(x, z) & = 1.
\end{align*}

Thus, condition (PCLI1) holds for any $x$ when $z < q$, and
the MP index is  $m(x) = x$ for $x < q$.

\paragraph{Case II: $q \leqslant z < h_{\infty}$.}
Consider the subcase $h_{t-1}(q) \leqslant   z < h_t(q)$ for $t \geqslant 1$.
The  metrics can be given in closed form in terms of 
$G_t(q) \triangleq G(q, h_{t-1}(q))$, $G_t(q + \rho)  \triangleq G(q + \rho, h_{t-1}(q))$, $F_t(q) \triangleq F(q, h_{t-1}(q))$ and $F_t(q + \rho)  \triangleq F(q + \rho, h_{t-1}(q))$, which are readily evaluated.
Thus, for $x > z$, 
\begin{align*}
F(x, z) & = x + \beta (x F_t(q + \rho) + (1-x) F_t(q)) \\
f(x, z) & = x - \beta (h(x) + \beta h(x) F_t(q + \rho) + \beta (1-h(x)) F_t(q) - x F_t(q + \rho) - (1-x) F_t(q)) \\
G(x, z) & = 1 + \beta (x G_t(q + \rho) + (1-x) G_t(q)) \\
g(x, z) & = 1 - \beta (1 + \beta h(x) G_t(q + \rho) + \beta (1-h(x)) G_t(q) - x G_t(q + \rho) - (1-x) G_t(q)).
\end{align*}

For $x \leqslant z$, letting $s \leqslant t+1$ be such that $h_{s-1}(x) \leqslant z < h_s(x)$, we have
\begin{align*}
F(x, z) & = \beta^s (h_s(x) + \beta h_s(x) F_t(q + \rho) + \beta \big(1-h_s(x)\big) F_t(q)) \\
f(x, z) & = x - \beta (\beta^{s-1}  [h_{s}(x) + \beta h_{s}(x)  F_t(q + \rho) 
+ \beta (1-h_{s}(x)) F_t(q)] -  x F_t(q + \rho) -  (1-x) F_t(q)) \\
G(x, z) & = \beta^s (1 + \beta h_s(x) G_t(q + \rho) + \beta \big(1-h_s(x)\big) G_t(q)) \\
g(x, z) & = 1 - \beta (\beta^{s-1}  [1 + \beta h_{s}(x)  G_t(q + \rho) 
+ \beta (1-h_{s}(x)) G_t(q)] -  x G_t(q + \rho) -  (1-x) G_t(q)). 
\end{align*}

Hence, for $q \leqslant x < h_{\infty}$, and provided $g(x, x) > 0$, the MP index is 
\[
m(x) =
\frac{x - \beta (h(x) + \beta h(x)  F_t(q + \rho) 
+ \beta (1-h(x)) F_t(q) -  x F_t(q + \rho) -  (1-x) F_t(q))}{1 - \beta (1 + \beta h(x)  G_t(q + \rho) 
+ \beta (1-h(x)) G_t(q) -  x G_t(q + \rho) -  (1-x) G_t(q))}.
\]

\paragraph{Case III: $h_{\infty} \leqslant z < q + \rho$.}
 In this 
case, the project metrics are given by
\begin{align*}
F(x, z) & = \frac{x}{1-\beta (q + \rho)} 1_{\{x > z\}} \\
f(x, z) & = \frac{x}{1-\beta (q + \rho)} 1_{\{h(x) \leqslant z\}} + \frac{(1-\beta  \rho) x -\beta  q}{1-\beta (q + \rho)} 
1_{\{h(x) > z\}} \\
G(x, z) & = \frac{1-\beta (q + \rho-x)}{1-\beta (q + \rho)} 1_{\{x > z\}} \\
g(x, z) & = 
\frac{1-\beta (q + \rho-x)}{1-\beta (q + \rho)} 1_{\{h(x) \leqslant z\}} +
\Big(1 -\beta +\beta  \frac{(1-\beta  \rho ) x -\beta  q}{1-\beta (q + \rho)}\Big) 1_{\{h(x) > z\}}.
\end{align*}

Hence, for $h_{\infty} \leqslant x < q + \rho$, as $h(x) \leqslant x$, the MP index is  
$
m(x) = x/(1-\beta (q + \rho-x)).
$

\paragraph{Case IV: $z \geqslant q + \rho$.} Finally, in this case we have
$F(x, z) = x 1_{\{x > z\}}$, $f(x, z) = x$, $G(x, z) = 1_{\{x > z\}}$, and $g(x, z)  = 
1$, 
and therefore, for $x \geqslant q + \rho$, the MP index is $m(x) = x$.

\begin{proposition}
\label{pro:pcloct}
The above optimal channel transmission model is PCL-indexable.
\end{proposition}
\proof 
Consider (PCLI1). 
The stronger result $g(x, z) \geqslant 1-\beta$ holds, which is trivial in cases I, III and IV.
In case II, for $h_{t-1}(q) < z < h_t(q)$, $g(\cdot, z)$ is piecewise linear \emph{c\`agl\`ad} (left-continuous with right limits) with $t+1$ increasing pieces, which follows from $G_t(q + \rho) >  G_t(q)$ and 
\[
\frac{\partial}{\partial x} g(x, z) = 
\begin{cases}
\beta  (1 - (\beta \rho)^{t+1}) (G_t(q + \rho) -  G_t(q)), & \textup{ for } 
0 < x < h_{-t}(z) \\
\beta  (1 - (\beta \rho)^{s}) (G_t(q + \rho) -  G_t(q)), & \textup{ for } 
h_{-s}(z) < x < h_{-s+1}(z), \, s = 2, \ldots, t \\
\beta  (1 - \beta \rho) (G_t(q + \rho) -  G_t(q)), & \textup{ for }  h_{-1}(z) < x < 1.
\end{cases}
\]
Hence, $g(x, z) > g(0^{\scriptscriptstyle +}, z)$ for $0 < x \leqslant h_{-t}(z)$, 
$g(x, z) > g(h_{-s}(z)^{\scriptscriptstyle +}, z)$ for $h_{-s}(z) < x \leqslant h_{-s+1}(z)$ with $2 \leqslant s \leqslant t$, and  $g(x, z) > g(h_{-1}(z)^{\scriptscriptstyle +}, z)$ for $h_{-1}(z) < x \leqslant 1$.
Now, it can be verified that $g(0^{\scriptscriptstyle +}, z) = 1$, and, using that 
$(\partial/\partial z) g(h_{-s}(z)^{\scriptscriptstyle +}, z) = 
\beta  (\rho^{-s} - \beta^s) (G_t(q + \rho) - G_t(q)) > 0$, $g(h_{-s}(z)^{\scriptscriptstyle +}, z) > g(h_{t-s-1}(q)^+, h_{t-1}(q))$ for $s = 1, \ldots, t$.
One can then show through algebraic and calculus arguments that $g(h_{t-s-1}(q)^+, h_{t-1}(q)) \geqslant 1-\beta$ for each such $s$, and thus establish satisfaction of (PCLI1).

As for condition (PCLI2), it is also easily verified in cases I, III and IV, viz., for $x \in [0, q) \cup [h_\infty, 1]$. 
Its verification in case II involves algebraic and calculus arguments.

Regarding (PCLI3), we use the sufficient conditions given in Section \ref{s:tcopwconstG}. Thus, taking 
\[
D(x) \triangleq \{h_t(x)\colon t \geqslant 0\} 
 \cup \{h_t(q)\colon t \geqslant 0\} \cup \{h_\infty, q, q + \rho\},
\]
Assumption 
\ref{ass:pwctraj} holds, whence Lemma \ref{lma:pwctraj} and Proposition \ref{pro:1cpcli3} imply (PCLI3).
\qquad
 \endproof

Thus, Theorem \ref{the:pcliii} gives that the model is threshold-indexable with Whittle index $m(\cdot)$.

\section{Concluding remarks.}
\label{s:conc}
This paper has presented PCL-based sufficient conditions for indexability along with optimality of threshold policies, for general restless bandits in a real-state discrete-time discounted setting,  demonstrating their validity. 
This extends the approach previously developed by the author for discrete-state bandits.
As for the value of the proposed approach relative to prevailing approaches based on first proving optimality of threshold policies (in a strong sense; see Section \ref{s:assub} and Section \ref{s:srsti}), the main example demonstrating it is the recent work in \citet{danceSi19}, which deploys the sufficient conditions presented here to  establish threshold-indexability of an important model that has not yielded to prevailing approaches, the Kalman filter target tracking bandit.   
Furthermore, the proposed approach applies to the larger class of threshold-indexable projects, whereas the conventional approach applies to the more restricted class of strongly threshold-indexable projects. 
We further stress that two  of the three sufficient conditions for threshold-indexability, viz.\ (PCLI2, PCLI3), are necessary (see Section \ref{s:potip}).  

Among the issues raised by this paper we highlight that further work is required to facilitate application of the proposed conditions to relevant model classes, by providing simpler means of checking their satisfaction than direct verification. Steps in that direction are already taken in this paper with condition (PCLI3). See Section \ref{s:ccpcli3}.  
It would also be desirable 
to identify properties of model primitives implying PCL-indexability. Extending the results to the average-reward criterion is another goal for further research. 


\bibliographystyle{plainnat}


\end{document}